\documentclass
[a4paper]{amsart}

\usepackage{amsthm,amsmath,amssymb}
\usepackage{extarrows}

\usepackage{tikz-cd}
\usepackage{quiver}

\numberwithin{equation}{section}

\newtheorem{corollary}[equation]{Corollary}

\newtheorem{lemma}[equation]{Lemma}

\newtheorem{proposition}[equation]{Proposition}

\newtheorem{theorem}[equation]{Theorem}
\theoremstyle{remark}
\newtheorem{remark}[equation]{Remark}
\theoremstyle{definition}
\newtheorem{construction}[equation]{Construction}
\newtheorem{convention}[equation]{Convention}
\newtheorem{definition}[equation]{Definition}
\newtheorem{example}[equation]{Example}

\newtheorem{notation}[equation]{Notation}
\newtheorem{terminology}[equation]{Terminology}

\newcommand{\Bir}{\mathbb{B}}

\newcommand{\boundary}{\partial}

\newcommand{\cat}{\mathcal}
\newcommand{\Cat}{\ccat{Cat}}

\newcommand{\Catset}{\Cat_1\times_{\Gpd}\Set}

\newcommand{\ccat}{\mathrm}
\newcommand*{\cimp}[1]{[\boldsymbol{#1}]}

\newcommand{\co}{\mathrm{co}}

\DeclareMathOperator*{\Compose}{\bigcirc}

\newcommand{\contra}{\mathrm{contra}}

\newcommand{\disj}{\sqcup}
\DeclareMathOperator*{\Disj}{\bigsqcup}

\newcommand{\Dual}{\mathbb{D}}

\DeclareMathOperator{\End}{End}
\renewcommand{\equiv}{\sim}

\newcommand{\equivwith}{\simeq}
\newcommand{\equivto}{\xrightarrow{\equiv}}

\newcommand{\from}{\leftarrow}
\newcommand{\Fun}{\ccat{Fun}}

\newcommand{\Gpd}{\ccat{Gpd}}

\newcommand{\Hom}{\mathrm{Hom}}

\newcommand{\id}{\mathrm{id}}
\newcommand{\contains}{\supset}
\newcommand{\content}{\underline}

\newcommand{\into}{\hookrightarrow}

\newcommand{\kore}{\textbf}
\DeclareMathOperator{\lad}{lad}

\newcommand{\Lattice}{\ccat{Lattice}}

\newcommand{\longequivto}{\xlongrightarrow{\equiv}}

\newcommand{\longto}{\longrightarrow}

\newcommand{\Mod}{\ccat{Mod}}

\newcommand{\noloc}{:\!}

\DeclareMathOperator{\Ob}{Ob}

\newcommand{\op}{\mathrm{op}}

\newcommand{\Ord}{\ccat{Ord}}
\newcommand*{\order}[1]{\Ord_\mathrm{#1}}
\newcommand{\Ordi}{\order{i}}
\newcommand{\Ordit}{\order{it}}
\newcommand{\Ordt}{\order{t}}

\newcommand{\Poset}{\ccat{Poset}}

\newcommand{\expression}{\mathsf}
\newcommand{\propP}{\expression{P}}

\DeclareMathOperator{\rad}{rad}

\newcommand*{\resto}[1]{|_{#1}}

\newcommand{\Set}{\ccat{Set}}

\newcommand{\tensor}{\otimes}
\newcommand{\Tensor}{\bigotimes}

\newcommand{\truth}{\mathrm{T}}

\newcommand{\unity}{\boldsymbol{1}}
\newcommand{\Univ}{\mathbf{Univ}}
\newcommand{\varv}{\expression{v}}

\title[On some categories of finite totally ordered sets]{%
  Notes on some categories related to that of finite totally ordered
  sets.}

\author{Matsuoka, Takuo}
\email{motogeometop@gmail.com}



\begin{document} 
\begin{abstract}
  The purpose of these notes is to collect in one place some facts on
  the category of finite totally ordered sets and some related
  categories.
  More specifically, we collect some results on them which will be
  useful for the study of iteratedly meta theories of algebra in the
  style of our work \cite{riron}, which is a kind of higher order
  universal algebra.
\end{abstract}

\maketitle

\tableofcontents

\section{Preliminaries}
\label{sec:category}

In this section, we cover some parts of category theory to prepare for
our exposition.

\subsection{Categories and groupoids}
\label{sec:category-groupoid}

In this section, we collect some basic definitions and facts.

\begin{definition}
  Given a category $\cat{C}$, its \kore{underlying} groupoid
  $\cat{C}_0$ is the
  groupoid consisting of the same objects as those of $\cat{C}$ and
  exactly the isomorphisms in $\cat{C}$ as morphisms, such that the
  composition in $\cat{C}_0$ is restricted from $\cat{C}$.
\end{definition}

\begin{remark}
  \label{rem:underlying-groupoid-counit}
  We have a functor $\cat{C}_0\to \cat{C}$ which sends each object and
  morphism in $\cat{C}_0$ to itself in $\cat{C}$.
\end{remark}

\begin{lemma}\label{lemma:subcategory}
  Let $\cat{C}$, $\cat{D}$ be categories, and $F\colon\cat{C}\to
  \cat{D}$ a functor.
  \begin{enumerate}
  \item If $F$ is faithful, then, for every category $\cat{X}$, the
    induced functor
    \[
      \Fun(\cat{X}, F)\colon \Fun(\cat{X}, \cat{C})\longto
      \Fun(\cat{X}, \cat{D})
    \]
    on the functor categories is faithful.
  \item If moreover the underlying functor $F_0\colon\cat{C}_0\to
    \cat{D}_0$ of groupoids is fully faithful, then $\Fun(\cat{X},
    F)_0$ on the underlying groupoids is fully faithful.
  \end{enumerate}
\end{lemma}

\begin{terminology}
  A functor $F$ is \kore{monic} if it satisfies both of the
  assumptions of Lemma.
\end{terminology}

\begin{remark}
  Such $F$ is said to be \emph{pseudomonic} in \cite{modulate}.
  Our terminology is for simplicity, since we do not see any other
  notion of monicity which may help in these notes.
\end{remark}

\begin{example}
  For a category $\cat{C}$, the canonical functor $\cat{C}_0\to
  \cat{C}$ of Remark~\ref{rem:underlying-groupoid-counit} is monic.
\end{example}

\begin{corollary}
  Let $F\colon\cat{C}\to \cat{D}$ be a monic functor.
  Then, for every category $\cat{X}$, the induced functor
  \[
    \Fun(\cat{X}, \cat{C})\longto
    \Fun(\cat{X}, \cat{D})\times_{\Fun(\cat{X}_0, \cat{D}_0)}
    \Fun(\cat{X}_0, \cat{C}_0)
  \]
  is monic, where the target is the fibre product in the
  bicategory of categories.
\end{corollary}
\begin{proof}
  It suffices to note that the assumption also implies that the
  projection
    \[
    \Fun(\cat{X}, \cat{D})\times_{\Fun(\cat{X}_0, \cat{D}_0)}
    \Fun(\cat{X}_0, \cat{C}_0)\longto \Fun(\cat{X}, \cat{D})
  \]
  is fully faithful.
\end{proof}

\begin{remark}
  Lemma~\ref{lemma:subcategory} has the converse that a functor
  $F\colon \cat{C}\to \cat{D}$ is monic if, for every small
  category $\cat{X}$, the functor $\Fun(\cat{X}, F)_0$ is fully
  faithful.
  In fact, just two samples of $\cat{X}$ suffice, e.g., the
  contractible groupoid and $\cimp{1}$ of
  Notation~\ref{nota:1-simplex}.
  (See Convention~\ref{conv:poset} or
  Example~\ref{ex:poset-of--1-types} on the latter sample.)
\end{remark}

\begin{definition}
  The \kore{underlying category} of a bicategory $\cat{C}$, is the
  category $\cat{C}_1$ enriched in the Cartesian bicategory of
  groupoids whose object is an object of $\cat{C}$, and such that
  \[
    \Hom_{\cat{C}_1}(X, Y) = \Hom_\cat{C}(X, Y)_0
  \]
  for every pair $X$, $Y$ of objects, and the composition in
  $\cat{C}_1$ is restricted from $\cat{C}$.
\end{definition}

\begin{remark}
  We have a canonical functor $\cat{C}_1\to \cat{C}$.
\end{remark}

\subsection{Partial orders}
\label{sec:poset}

\begin{notation}
  $\Cat$ denotes the bicategory of small categories (see Remark to
  follow).
\end{notation}

\begin{remark}
  By \kore{smallness}, we always mean essential smallness.
  In fact, \textbf{we never try to use or introduce any notion which
    is not invariant under equivalence of categories}.
\end{remark}

\begin{notation}
  $\Gpd$ denotes the full sub-bicategory of $\Cat$ consisting of
  small groupoids.
\end{notation}

$\Gpd$ is in fact in $\Cat_1$.
The functor $(-)_0\colon \Cat_1\to \Gpd$ is the right adjoint of the
inclusion with counit given by the canonical functor $\cat{C}_0\to
\cat{C}$ of Remark~\ref{rem:underlying-groupoid-counit} for each
$\cat{C}\in\Cat_1$.

\begin{notation}
  $\Set$ denotes the category of sets.
\end{notation}

Any category can be regarded as a bicategory by considering its hom
sets as the corresponding homotopically discrete groupoids.
$\Set$ regarded as a bicategory in this manner embeds in $\Gpd$ as
the full subcategory consisting of small homotopically discrete
groupoids.
Indeed, a small homotopically discrete groupoid can be identified with
the set of its isomorphism classes.

In this section, we recall extention of this equivalence to one
between partially ordered sets and some small categories.

\begin{terminology}
  A \kore{homotopy ($-1$)-type} is a set in which all elements are
  equal.
\end{terminology}

Homotopy ($-1$)-types are handy truth values to come from the relation
between logic and sets as follows.
\begin{definition}
  \label{def:truth-value}
  Given a proposition $\propP$ which we assume for simplicity
  can be expressed without using the variable ``$t$'', let
  \begin{equation}
    \label{eq:truth-value}
    \truth(\propP)\;=\;\{t \in T\: \mid\: \content{\propP}\},
  \end{equation}
  where $T$ is a terminal set, and the expression
  ``$\content{\propP}$'' does not specify an object, but should be
  replaced by the content of $\propP$ expressed without using the
  variable ``$t$''.
  (For instance, if $\propP$ is the proposition ``$T=0$'' in a context
  where this expression makes sense, then $\truth(\propP)=\{t\in
  U\:\mid\:T=0\}$ in the context, where $U$ is a terminal set.)
\end{definition}

\begin{remark}
  Choice of the terminal set $T$ in Definition is irrelevant since the
  right hand side of (\ref{eq:truth-value}) is functorial in $T$.
\end{remark}

$\truth(\propP)$ is a homotopy ($-1$)-type, and indicates whether
$\propP$ is true.
More generally, homotopy ($-1$)-types are the possible fibres of the
inclusion map of a subset of any set, which corresponds to a predicate
on an element of the set.

\begin{remark}
  To extend the definition of $\truth(\propP)$ for any proposition
  $\propP$, let $\varv$ be a variable (which we assumed to be ``$t$''
  in Definition~\ref{def:truth-value}) not necessary for expressing
  $\propP$.
  Then $\truth(\propP)
  = \{\content{\varv}\in T\: \mid\: \content{\propP}\}$, where the
  expressions ``$\content{\varv}$'' and ``$\content{\propP}$'' should
  respectively be replaced by the symbol of $\varv$ and an expression
  for $\propP$ not containing the variable $\varv$.
\end{remark}

\begin{construction}\label{constr:category-from-poset}
  From a partially ordered set $P$, we obtain a category $\cat{P}$ whose object is an
  element of $P$, and such that
  \[
    \Hom_\cat{P}(x, y)\;=\;\truth(\propP(x, y))
  \]
  for every pair $x$, $y$ of objects of $\cat{P}$, where, for elements
  $x$, $y$ of $P$, $\propP(x, y)$ denotes the proposition
  \[
    x \le y\quad\text{in $P$}.
  \]
\end{construction}
A small category which is equivalent to the one obtained from a partially ordered set
in this manner is characterized by the property that the set of
morphisms in it is a homotopy ($-1$)-type for every pair of source and
target.

\begin{notation}
  $\Poset$ denotes the full sub-bicategory of $\Cat$ consisting of
  these small categories.
\end{notation}

\begin{terminology}
  The \kore{$2$-category of partially ordered sets} is the bicategory
  $\Poset$.
\end{terminology}

The groupoid of morphisms between every pair of source and target in
$\Poset_1$ is homotopically discrete, so $\Poset_1$ can be considered
as a category in the usual, unenriched sense, i.e., one enriched in
$\Set$.
Moreover, between the objects of $\Poset_1$ constructed from partially ordered sets,
the homotopically discrete groupoid of morphisms in $\Poset_1$ can be
identified with the set of non-increasing maps between those partially ordered sets.
Thus, Construction~\ref{constr:category-from-poset} is an
equivalence to $\Poset_1$ from the usual category of partially ordered sets (with
non-increasing maps as morphisms in it).

\begin{convention}
  \label{conv:poset}
  \textbf{We identify the category of partially ordered sets with
    $\Poset_1$ through this equivalence, and identify every obeject of
    the latter with the corresponding partially ordered set}, so, by a ``poset'', we
  normally mean an object of $\Poset$.
\end{convention}

\begin{example}
  \label{ex:poset-of--1-types}
  The category of homotopy ($-1$)-types may be identified with the
  partially ordered set $\cimp{1}$ of Notation~\ref{nota:1-simplex}.
\end{example}

In fact, the category of morphisms between every pair of source and
target in $\Poset$ belongs to $\Poset$, so $\Poset$ can be considered
as a category enriched in the Cartesian category $\Poset_1$.
This is the enrichment of $\Poset_1$ in itself (in
Terminology~\ref{term:enrich}) coming from its Cartesian closedness
(see Remark~\ref{rem:canonical-forget}).
Note that the underlying groupoid of a poset is homotopically
discrete, and is the underlying set of the partially ordered set in the usual sense.

\subsection{Enrichment}
\label{sec:enrich}

In this section, we introduce some terms.

\begin{terminology}
  \label{term:enrich}
  Let $F\colon \cat{V}\to \cat{W}$ be a lax monoidal functor between
  monoidal categories.

  When $F$ is considered a forgetful functor, for a category $\cat{C}$
  enriched in $\cat{V}$, the $\cat{W}$-enriched category
  \kore{underlying} $\cat{C}$, with the lax monoidal structure of $F$
  understood, is the $\cat{W}$-enriched category induced from
  $\cat{C}$ through $F$.

  Conversely, given a category $\cat{D}$ enriched in $\cat{W}$, an
  \kore{enrichment} of $\cat{D}$ in $\cat{V}$ is a $\cat{V}$-enriched
  category which $\cat{D}$ underlies.
\end{terminology}

\begin{remark}
  \label{rem:canonical-forget}
  Often, e.g., in \cite{enrich}, $F$ is the functor
  $\cat{V}\to \Set$ represented by the unit object $\unity$ of
  $\cat{V}^\op$, where the lax monoidal structure comes from the
  monoid structure of $\unity$.
  In particular, a category enriched in $\cat{W}$ is a category in the
  usual sense.
  This forgetful functor is one with respect to which a closed
  symmetric monoidal category is enriched in itself.

  For the general case, it may often be reasonable to require $F\colon
  \cat{V}\to \cat{W}$ to commute with the respective forgetful
  functors to $\Set$ on $\cat{V}$ and $\cat{W}$.
\end{remark}

\begin{terminology}
  The \kore{canonical forgetful functor} on a monoidal category
  $\cat{V}$ is the lax monoidal functor to $\Set$ represented by the
  unit monoid $\unity$ in $\cat{V}^\op$.
\end{terminology}

This functor may \emph{not} be useful on every monoidal category
outside the context of enrichment.

\begin{example}
  The canonical forgetful functor on a co-Cartesian monoidal category
  is constant at the terminal monoid.  That is, it factorizes through
  the lax monoidal functor from the terminal monoidal category to
  $\Set$ classifying the terminal monoid in $\Set$.
\end{example}

Validity of the following claim is straightforward to establish.

\begin{proposition}
  Let $F\colon \cat{V}\to \cat{W}$ be a lax monoidal functor of
  monoidal cagegories commuting with the canonical forgetful functors.
  Let $f\colon \cat{C}\to \cat{D}$ be a functor of $\cat{V}$-enriched
  categories such that the $\cat{W}$-enriched functor $F_*f\colon
  F_*\cat{C}\to F_*\cat{D}$ admits a left (say) adjoint.
  Then $f$ admits a left adjoint such that the adjunction induces the
  $\cat{W}$-enriched adjunction \textbf{if} $F$ reflects isomorphisms.
\end{proposition}

\section{Birkhoff duality}
\label{sec:duality}

In this section, we study the Birkhoff duality for finite totally
ordered sets as an equivalence of categories.

\subsection{The duality}
\begin{notation}
  $\Ord$ denotes the full subcategory of $\Poset_1$ consisting of
  finite totally ordered sets.
\end{notation}

\begin{remark}
  The category $\Ord$ is useful for universal algebra since it is the
  monoidal category where the universal monoid is.
\end{remark}

The monoidal structure will be looked at in
Section~\ref{sec:operations}.
In this section, we shall only consider its structure of a category.

The underlying groupoid $\Ord_0$ has discrete homotopy type, which may
be described e.g., as the set of natural numbers.

Let us consider another category, which happens to have monic
forgetful functor to $\Ord$.

\begin{definition}\label{def:ordit}
  The category $\Ordit$ is specified as follows.
  \begin{itemize}
  \item An object of $\Ordit$ is a \emph{non-empty} totally ordered
    set.
  \item A morphism in it between such objects is a non-increasing map
    between them which preserves both the minimal element and the
    maximal element.
  \end{itemize}
\end{definition}

The forgetful functor $\Ordit\to\Ord$ is monic.
The fully faithful functor induced on the underlying groupoids can be
identified with an injective map of sets, which can be understood
e.g., as the inclusion of strictly positive integers into all natural
numbers.

As a poset, each object of $\Ordit$ is a lattice, and $\Ordit$ is in
fact a full subcategory of the category of bounded lattices.
Moreover, a very non-Boolean part of the Birkhoff duality (between
finite distributive lattices and finite posets) gives an equivalence
$\Ordit^\op\equivwith\Ord$.
(The Boolean part of the duality is an equivalence of finite
Boolean algebras with finite sets, while objects of $\Ord$ are in a
way furthest from sets among finite posets.)

For the rest of this section, we recall the equivalence, which is
given ``by a dualizing object''.

\begin{notation}\label{nota:1-simplex}
  $\cimp{1}$ denotes the initial object of $\Ordit$.
\end{notation}

We would like to give:
\begin{itemize}
\item the object $\cimp{1}$ of $\Ordit$ the structure of an
  $\Ord$-object relative to the forgetful functor $\Ord\to\Set$, which
  forgets the order and the finiteness, and
\item the object $\cimp{1}$ of $\Ord$ ($\contains \Ordit$) the
  structure of an $\Ordit$-object relative to the composite
  $\Ordit\to\Ord\to\Set$ of fogetful functors.
\end{itemize}

To prepare for this, let us enrich $\Ord$ in $\Poset_1$ with respect
to the lax monoidal forgetful functor $\Poset_1\to \Set$.
We do this by letting (the enriched version of) $\Ord$ be a full
$\Poset_1$-enriched subcategory of $\Poset_1$.

\begin{notation}
  For objects $I$, $J$ of $\Ord$, $\Hom_\Ord(I, J)$ denotes the hom
  poset in this enrichment of $\Ord$.
\end{notation}

In other words, $\Hom_\Ord(I, J)$ denotes $\Hom_{\Poset}(I, J)$
considered as a poset.
The underlying set of $\Hom_\Ord(I, J)$ is as desired the \emph{set}
of morphisms in $\Ord$ from which we started.

Recall that we wanted to lift the presheaf on $\Ord$ represented by
$\cimp{1}$ through forgetful functors.
Note that the forgetful functor $\Ord\to \Set$ factorizes as the
inclusion $\Ord\into\Poset_1$ followed by the forgetful functor
$\Poset_1\to \Set$.
We have in particular lifted the presheaf to the functor $\Hom_\Ord(-,
\cimp{1})\colon \Ord^\op\to \Poset_1$ against the forgetful functor
$\Poset_1\to\Set$, where the functoriality of $\Hom_\Ord(-, \cimp{1})$
results from the closedness of $\Poset_1$.
(The functoriality is also the one obtained by using the
representation of the lax monoidal forgetful functor $\Poset_1\to
\Set$ by the unit monoid in $\Poset_1^\op$.)
Moreover, this lift lands in the full subcategory $\Ord$ of
$\Poset_1$.
Thus, we have lifted the presheaf to a contravariant functor to
$\Ord$.
Finally, this contravariant functor comes from a functor
$\Ord^\op\to \Ordit$.
(Note Lemma~\ref{lemma:subcategory} or its corollary.)
This achieves the second of the objectives we posed.

To work for the other objective, in order to lift the presheaf on
$\Ordit$ represented by $\cimp{1}$ to a contravariant functor to
$\Poset_1$, it suffices to enrich $\Ordit$ in $\Poset_1$.
We do this by enriching it in such a way that the forgetful functor
$\Ordit\to \Ord$ is enriched in such a manner that, for every pair
$I$, $J$ of objects of $\Ordit$, the induced map
\[
  \Hom_{\Ordit}(I, J)\longto \Hom_\Ord(I, J)
\]
of the hom posets will be fully faithful.
(Recall that a map of posets is a functor for us.)
This is possible since what we need to make the underlying functor of
the desired fully faithful functor is indeed an injection, i.e., a
fully faithful functor of homotopically discrete groupoids.
Moreover, the resulting lift $\Hom_{\Ordit}(-,
\cimp{1})\colon\Ordit^\op\to \Poset_1$ lands in the full subcategory
$\Ord$ of $\Poset_1$.
This achieves the first of our objectives.

Thus, between $\Ord$ and $\Ordit$, we have obtained a contravariant
functor in each direction which lifts a representable presheaf through
the forgetful functor.
Moreover, the representing objects of the presheaves have common
underlying set.

To see that these functors are mutually inverse equivalences given by
the Birkhoff duality, since the enrichment in $\Poset_1$ which we have
given to $\Ordit$
is a restriction of the usual enrichment of the category $\Lattice$ of
bounded lattices, it suffices to note that the forgetful functor
$\Lattice\to \Poset_1$ is also monic.

\subsection{Intermediate self duality}
\label{sec:self-duality}

In this section, we put the Birkhoff duality in a more structured
duality.

\begin{notation}
  $\Ordt$ denotes the category of non-empty totally ordered sets and
  non-increasing maps which preserves the maximal element, defined
  similarly to the way how $\Ordit$ was defined in
  Definition~\ref{def:ordit}.
\end{notation}

\begin{remark}
  \label{rem:universal-monad-algebra}
  $\Ordt$ is relevant to universal algebra since it is a module
  category over $\Ord$ where the universal module over the universal
  monoid is.
  (What determines the side of the module is the conventions leading
  to the monoidal variance of the module structure functor
  $\Ord\to \End(\Ordt)$ together with the convention on on which side
  $\Ordt$ is a module over $\End(\Ordt)$ (resulting in arbitrary
  attachment of the names of sides to the variances).
  However, $\Ord$ is equivalent to its monoidal opposite.
  A side/variance-free description of a (left or right) $\Ord$-module
  structure on a category $\cat{C}$ is that it is a monad on
  $\cat{C}$, and a module over the universal monoid in an
  $\Ord$-module category $\cat{C}$ will be an algebra over the
  structure monad on $\cat{C}$.)
\end{remark}

The structure of an $\Ord$-module on $\Ordt$ will be looked at in
Section~\ref{sec:ord-module}.

The forgetful functor $\Ordt\to \Ord$ is monic.
It follows that the forgetful functor $\Ordit\to \Ordt$ is monic.

$\cimp{1}$ will be a dualizing object in $\Ordt$ as well.
To see this, enrich $\Ordt$ in $\Poset_1$ similarly as we did
$\Ordit$.
Then the functor $\Hom_{\Ordt}(-, \cimp{1})\colon \Ordt^\op\to
\Poset_1$ to lift the presheaf on $\Ordt$ represented by $\cimp{1}$
lands in $\Ordt$.

\begin{notation}
  $\Dual\colon \Ordt^\op\to\Ordt$ denotes the resulting functor.
\end{notation}

\begin{theorem}
  $\Dual$ is inverse to $\Dual\colon \Ordt\to\Ordt^\op$.
\end{theorem}
\begin{proof}
  These functors are adjoint to each other with unit map
  $\epsilon\colon \id\to \Dual\Dual$ coinciding with the counit map.
  Therefore, it suffices to prove that $\epsilon$ is an isomorphism.
  The map $\epsilon\Dual\colon \Dual\to \Dual\Dual\Dual$ is a section,
  but for every object $I\in \Ob\Ordt$, we have $\Dual I\equivwith I$.
  The result follows immediately.
\end{proof}

We shall relate this duality to the Birkhoff duality.

\begin{notation}
  $\Bir\colon \Ord^\op\rightleftarrows\Ordit\noloc \Bir$ denote the
  Birkhoff duality functors.
  (We use the same symbol for the two functors.
  This is an abbreviation for writing one as $\Bir$ and the other as
  $\Bir^{-1}$.)
\end{notation}

Note that each of the forgetful functors $\Ordit\into \Ordt\into \Ord$
admits a left adjoint.

\begin{notation}
  The adjoints are denoted as in $\Ordit\xleftarrow{i}
  \Ordt\xleftarrow{t} \Ord$.
\end{notation}

The adjunctions are in fact enriched in $\Poset_1$.

Note that the categories $\Ord$, $\Ordt$, $\Ordit$ lie in an
essentially \emph{$\Set$-enriched}
category, e.g., $\Catset$.

\begin{proposition}\label{prop:sequence-duality}
The squares
% https://q.uiver.app/#q=WzAsOSxbMiwxLCJcXE9yZGl0Xlxcb3AiXSxbMSwxLCJcXE9yZHReXFxvcCJdLFsyLDIsIlxcT3JkIl0sWzEsMiwiXFxPcmR0Il0sWzAsMSwiXFxPcmReXFxvcCJdLFswLDIsIlxcT3JkaXQiXSxbMCwwLCJcXE9yZGl0Il0sWzEsMCwiXFxPcmR0Il0sWzIsMCwiXFxPcmQiXSxbMCwyLCJcXEJpciJdLFsxLDMsIlxcRHVhbCJdLFsxLDAsImkiXSxbMywyLCIiLDEseyJzdHlsZSI6eyJ0YWlsIjp7Im5hbWUiOiJob29rIiwic2lkZSI6InRvcCJ9fX1dLFs0LDUsIlxcQmlyIl0sWzQsMSwidCJdLFs1LDMsIiIsMix7InN0eWxlIjp7InRhaWwiOnsibmFtZSI6Imhvb2siLCJzaWRlIjoidG9wIn19fV0sWzYsNCwiXFxCaXIiXSxbNywxLCJcXER1YWwiXSxbOCwwLCJcXEJpciJdLFs2LDcsIiIsMSx7InN0eWxlIjp7InRhaWwiOnsibmFtZSI6Imhvb2siLCJzaWRlIjoidG9wIn19fV0sWzcsOCwiIiwxLHsic3R5bGUiOnsidGFpbCI6eyJuYW1lIjoiaG9vayIsInNpZGUiOiJ0b3AifX19XV0=
\[\begin{tikzcd}
	\Ordit & \Ordt & \Ord \\
	{\Ord^\op} & {\Ordt^\op} & {\Ordit^\op} \\
	\Ordit & \Ordt & \Ord
	\arrow[hook, from=1-1, to=1-2]
	\arrow["\Bir", from=1-1, to=2-1]
	\arrow[hook, from=1-2, to=1-3]
	\arrow["\Dual", from=1-2, to=2-2]
	\arrow["\Bir", from=1-3, to=2-3]
	\arrow["t", from=2-1, to=2-2]
	\arrow["\Bir", from=2-1, to=3-1]
	\arrow["i", from=2-2, to=2-3]
	\arrow["\Dual", from=2-2, to=3-2]
	\arrow["\Bir", from=2-3, to=3-3]
	\arrow[hook, from=3-1, to=3-2]
	\arrow[hook, from=3-2, to=3-3]
\end{tikzcd}\]
commute, where the vertical arrows are respective duality
isomorphisms.
\end{proposition}
\begin{proof}
  We show the commutativity of the lower squares.
  This suffices since each of the vertical arrows in the upper squares
  is the inverse of the functor right below it.

  To show that the left square commutes, it suffices to show that the
  two functors are equalized by the monic forgetful functor
  $\Ordt\into \Poset_1$.
  However, the two composites with the forgetful functor are
  isomorphic by the $\Poset_1$-enriched adjunction isomorphism.

  The right square then commutes since it is the square formed by the
  adjoints of the functors in the left square.
  Alternatively, the square can directly be shown to commute by
  similar arguments as we have given for the left square.
\end{proof}

Thus, the dualities give an isomorphism of the sequence
$\Ord\xrightarrow{t} \Ordt\xrightarrow{i} \Ordit$ with
$\Ordit^\op\into \Ordt^\op\into \Ord^\op$.

In the left squares of Proposition, each of the Birkhoff duality
functor $\Bir$ post-composed with a monic functor is described in a
way which does not rely on either of the functors.
This in fact leads to a description not relying on $\Bir$, of $\Bir$
itself (or equivalently, of $\Bir$ composed with an unrelated
isomorphism).
Let us find the description.
Essentially, we would like to describe the image of the monic functor
with which we post-composed with $\Bir$ previously without using the
monic functor or its source.

\begin{notation}
  $\Ordi$ denotes the category of non-empty totally ordered sets and
  non-increasing maps which preserves the minimal element.
\end{notation}

Even though this category is equivalent to $\Ordt$ by the functor to
take the opposite of posets (and categories), we consider it as a
different (non-full) \emph{sub}category of $\Ord$ which gets
exchanged with $\Ordt$ by the said functor on $\Ord$.

\begin{notation}\label{nota:op}
  $\op$ denotes the functor to take the opposites at each place in the
  diagram
  % https://q.uiver.app/#q=WzAsOCxbMCwxLCJcXE9yZGl0Il0sWzEsMCwiXFxPcmR0Il0sWzEsMiwiXFxPcmRpIl0sWzIsMSwiXFxPcmQiXSxbMSwxLCJcXE9yZGl0Il0sWzMsMSwiXFxPcmQiXSxbMiwwLCJcXE9yZGkiXSxbMiwyLCJcXE9yZHQiXSxbMCw0LCJcXG9wIl0sWzMsNSwiXFxvcCJdLFsxLDYsIlxcb3AiXSxbMiw3LCJcXG9wIiwyXSxbMCwxLCIiLDAseyJzdHlsZSI6eyJ0YWlsIjp7Im5hbWUiOiJob29rIiwic2lkZSI6InRvcCJ9fX1dLFsxLDMsIiIsMix7ImxhYmVsX3Bvc2l0aW9uIjo3MCwic3R5bGUiOnsidGFpbCI6eyJuYW1lIjoiaG9vayIsInNpZGUiOiJ0b3AifX19XSxbMCwyLCIiLDIseyJzdHlsZSI6eyJ0YWlsIjp7Im5hbWUiOiJob29rIiwic2lkZSI6InRvcCJ9fX1dLFsyLDMsIiIsMCx7ImxhYmVsX3Bvc2l0aW9uIjo3MCwic3R5bGUiOnsidGFpbCI6eyJuYW1lIjoiaG9vayIsInNpZGUiOiJ0b3AifX19XSxbNCw2LCIiLDIseyJzdHlsZSI6eyJ0YWlsIjp7Im5hbWUiOiJob29rIiwic2lkZSI6InRvcCJ9fX1dLFs0LDcsIiIsMix7InN0eWxlIjp7InRhaWwiOnsibmFtZSI6Imhvb2siLCJzaWRlIjoidG9wIn19fV0sWzcsNSwiIiwyLHsic3R5bGUiOnsidGFpbCI6eyJuYW1lIjoiaG9vayIsInNpZGUiOiJ0b3AifX19XSxbNiw1LCIiLDIseyJzdHlsZSI6eyJ0YWlsIjp7Im5hbWUiOiJob29rIiwic2lkZSI6InRvcCJ9fX1dXQ==
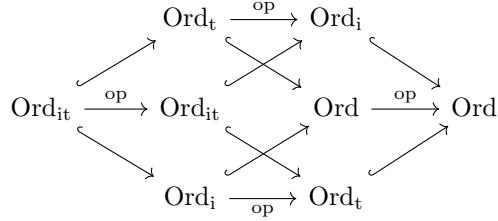
\[\begin{tikzcd}
	& \Ordt & \Ordi \\
	\Ordit & \Ordit & \Ord & \Ord \\
	& \Ordi & \Ordt
	\arrow["\op", from=1-2, to=1-3]
	\arrow[hook, from=1-2, to=2-3]
	\arrow[hook, from=1-3, to=2-4]
	\arrow[hook, from=2-1, to=1-2]
	\arrow["\op", from=2-1, to=2-2]
	\arrow[hook, from=2-1, to=3-2]
	\arrow[hook, from=2-2, to=1-3]
	\arrow[hook, from=2-2, to=3-3]
	\arrow["\op", from=2-3, to=2-4]
	\arrow[hook, from=3-2, to=2-3]
	\arrow["\op"', from=3-2, to=3-3]
	\arrow[hook, from=3-3, to=2-4]
\end{tikzcd}\]
\end{notation}

\begin{notation}
  The left adjoints of the forgetful functors $\Ordit\into \Ordi\into
  \Ord$ are denoted as in $\Ordit\xleftarrow{t} \Ordi\xleftarrow{i}
  \Ord$.
\end{notation}

\begin{notation}
  $\Dual\colon \Ordi^\op\to \Ordi$ denotes the duality isomorphism.
\end{notation}

\begin{lemma}\label{lem:functor-to-pull-back}
  The diagram
  % https://q.uiver.app/#q=WzAsNixbMCwxLCJcXE9yZF5cXG9wIl0sWzEsMiwiXFxPcmRpXlxcb3AiXSxbMiwyLCJcXE9yZGkiXSxbMywxLCJcXE9yZCJdLFsxLDAsIlxcT3JkdF5cXG9wIl0sWzIsMCwiXFxPcmR0Il0sWzAsMSwiaSIsMl0sWzEsMiwiXFxEdWFsIiwyXSxbMiwzLCIiLDIseyJzdHlsZSI6eyJ0YWlsIjp7Im5hbWUiOiJob29rIiwic2lkZSI6InRvcCJ9fX1dLFs1LDMsIiIsMCx7InN0eWxlIjp7InRhaWwiOnsibmFtZSI6Imhvb2siLCJzaWRlIjoidG9wIn19fV0sWzQsNSwiXFxEdWFsIl0sWzAsNCwidCJdXQ==
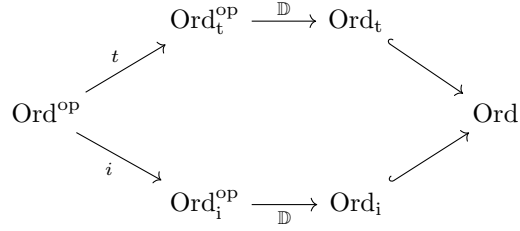
\[\begin{tikzcd}
	& {\Ordt^\op} & \Ordt \\
	{\Ord^\op} &&& \Ord \\
	& {\Ordi^\op} & \Ordi
	\arrow["\Dual", from=1-2, to=1-3]
	\arrow[hook, from=1-3, to=2-4]
	\arrow["t", from=2-1, to=1-2]
	\arrow["i"', from=2-1, to=3-2]
	\arrow["\Dual"', from=3-2, to=3-3]
	\arrow[hook, from=3-3, to=2-4]
\end{tikzcd}\]
commutes.
\end{lemma}
\begin{proof}
  It follows from Proposition~\ref{prop:sequence-duality} that the
  upper composite is equal to the composite
  % https://q.uiver.app/#q=WzAsMyxbMCwwLCJcXE9yZF5cXG9wIl0sWzEsMCwiXFxPcmRpdCJdLFsyLDAsIlxcT3JkIl0sWzEsMiwiIiwwLHsic3R5bGUiOnsidGFpbCI6eyJuYW1lIjoiaG9vayIsInNpZGUiOiJ0b3AifX19XSxbMCwxLCJcXEJpciJdXQ==
\[\begin{tikzcd}
	{\Ord^\op} & \Ordit & \Ord
	\arrow["\Bir", from=1-1, to=1-2]
	\arrow[hook, from=1-2, to=1-3]
\end{tikzcd}\]
to which the lower composite is also equal by the version of the same
proposition for $\Ordi$.
\end{proof}

\begin{lemma}\label{lem:pull-back}
  The square
  % https://q.uiver.app/#q=WzAsNCxbMCwxLCJcXE9yZGl0Il0sWzEsMCwiXFxPcmR0Il0sWzIsMSwiXFxPcmQiXSxbMSwyLCJcXE9yZGkiXSxbMSwyLCIiLDAseyJzdHlsZSI6eyJ0YWlsIjp7Im5hbWUiOiJob29rIiwic2lkZSI6InRvcCJ9fX1dLFszLDIsIiIsMix7InN0eWxlIjp7InRhaWwiOnsibmFtZSI6Imhvb2siLCJzaWRlIjoidG9wIn19fV0sWzAsMSwiIiwwLHsic3R5bGUiOnsidGFpbCI6eyJuYW1lIjoiaG9vayIsInNpZGUiOiJ0b3AifX19XSxbMCwzLCIiLDIseyJzdHlsZSI6eyJ0YWlsIjp7Im5hbWUiOiJob29rIiwic2lkZSI6InRvcCJ9fX1dXQ==
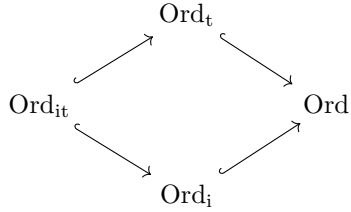
\[\begin{tikzcd}
	& \Ordt \\
	\Ordit && \Ord \\
	& \Ordi
	\arrow[hook, from=1-2, to=2-3]
	\arrow[hook, from=2-1, to=1-2]
	\arrow[hook, from=2-1, to=3-2]
	\arrow[hook, from=3-2, to=2-3]
\end{tikzcd}\]
is a pull-back square.
\end{lemma}
\begin{proof}
  The underlying groupoid of the fibre product is the homotopically
  discrete groupoid corresponding to the fibre product of the
  homotopically discrete underlying groupoids taken in the category
  of sets.
  The set of maps in the fibre product is the fibre product of the
  sets of maps.
\end{proof}

It follows from Proposition~\ref{prop:sequence-duality} that
$\Bir\colon \Ord^\op\to\Ordit$ is the functor to the fibre product
which corresponds to the commutative diagram of
Lemma~\ref{lem:functor-to-pull-back}.

We can get a similar description of $\Bir\colon \Ordit^\op\to \Ord$,
but to give a clearer view, the duality isomorphism of sequences of
Proposition~\ref{prop:sequence-duality} and its counterpart for
$\Ordi$ glue together to an isomorphism of the pull-back square of
Lemma~\ref{lem:pull-back} with the square
\begin{equation}
  \label{eq:dual-pull-back}
  % https://q.uiver.app/#q=WzAsNCxbMSwwLCJcXE9yZHReXFxvcCJdLFsyLDEsIlxcT3JkaXReXFxvcCJdLFsxLDIsIlxcT3JkaV5cXG9wIl0sWzAsMSwiXFxPcmReXFxvcCJdLFswLDEsImkiXSxbMiwxLCJ0Il0sWzMsMCwidCJdLFszLDIsImkiXV0=
\begin{tikzcd}
	& {\Ordt^\op} \\
	{\Ord^\op} && {\Ordit^\op} \\
	& {\Ordi^\op}
	\arrow["i", from=1-2, to=2-3]
	\arrow["t", from=2-1, to=1-2]
	\arrow["i", from=2-1, to=3-2]
	\arrow["t", from=3-2, to=2-3]
\end{tikzcd}
\end{equation}

\begin{remark}
  As shown in Notation~\ref{nota:op}, the square of
  Lemma~\ref{lem:pull-back} is isomorphic to its flipped self by the
  isomorphisms $\op$.
  The isomorphisms $\op$ also flip the square
  (\ref{eq:dual-pull-back}).
  The isomorphism between the flipped squares obtained by conjugating
  the duality isomorphism by $\op$ is the same duality isomorphism
  (flipped).
  Indeed, the duality isomorphisms commute with $\op$ ($\op$ is
  ``dual'' to $\op$) since $\cimp{1}^\op=\cimp{1}$ (with nothing
  special of $\cimp{1}$).
  In the next section, we shall give a direct description of the
  composite of the duality isomorphism with $\op$.
\end{remark}

\subsection{Description via the adjoint maps}

In this section, we give another description of the dualities we have
considered above.

There is another equivalence
\begin{equation}
  \label{eq:equivalence-to-take-adjoints}
  \lad\colon\Ordt^\op\xlongleftrightarrow{\equiv} \Ordi\noloc\rad,
\end{equation}
where the functors are defined as follows.

\begin{definition}
  The functors (\ref{eq:equivalence-to-take-adjoints}) map every
  object to itself.
  For a morphism, i.e., an appropriate functor of posets, $\lad$
  (resp.~$\rad$) maps it to its left (resp.~right) adjoint.
\end{definition}

\begin{theorem}\label{thm:adjoint-dual}
  The diagram
  % https://q.uiver.app/#q=WzAsNixbNCwwLCJcXE9yZGleXFxvcCJdLFs1LDEsIlxcT3JkaSJdLFszLDEsIlxcT3JkdCJdLFsyLDAsIlxcT3JkdF5cXG9wIl0sWzAsMCwiXFxPcmRpXlxcb3AiXSxbMSwxLCJcXE9yZGkiXSxbMCwxLCJcXER1YWwiXSxbMiwxLCJcXG9wIiwyXSxbMywwLCJcXG9wIl0sWzMsMiwiXFxEdWFsIl0sWzAsMiwiXFxyYWQiXSxbNSwyLCJcXG9wIiwyXSxbNCwzLCJcXG9wIl0sWzQsNSwiXFxEdWFsIiwyXSxbMyw1LCJcXGxhZCIsMl1d
\[\begin{tikzcd}
	{\Ordi^\op} && {\Ordt^\op} && {\Ordi^\op} \\
	& \Ordi && \Ordt && \Ordi
	\arrow["\op", from=1-1, to=1-3]
	\arrow["\Dual"', from=1-1, to=2-2]
	\arrow["\op", from=1-3, to=1-5]
	\arrow["\lad"', from=1-3, to=2-2]
	\arrow["\Dual", from=1-3, to=2-4]
	\arrow["\rad", from=1-5, to=2-4]
	\arrow["\Dual", from=1-5, to=2-6]
	\arrow["\op"', from=2-2, to=2-4]
	\arrow["\op"', from=2-4, to=2-6]
\end{tikzcd}\]
commutes.
\end{theorem}
\begin{proof}
  Commutativity of the rightmost triangle suffices by symmetry (or the
  commutativity of any parallelograms could be used in addition).
  Thus, suppose given a morphism $f^*\colon I\to J$ in $\Ordi$.
  We would like to prove that the map $\Dual(f^*)\colon \Dual J\to
  \Dual I$ in $\Ordi$ coincides with the map $f_*^\op\colon J^\op\to
  I^\op$, where $f_*$ denotes the map $\rad(f^*)\colon J\to I$ in
  $\Ordt$.
  ($f_*^\op$ is the same as $f_*$ in $\Ord$, but we chose this
  notation hoping to make it clearer that the morphism is considered
  in $\Ordi$.)

  We may assume that $J = \cimp{1}$.
  Indeed, suppose we have the result in this case.
  Then, for the   general case, all we need to prove is that, for
  every morphism
  $j_*\colon \cimp{1}\to J$ in $\Ordt$, the maps $f_*^\op$ and
  $\Dual(f^*)$ are equalized by $j_*^\op\colon\cimp{1}^\op\to J^\op$.
  However, in the diagram
  % https://q.uiver.app/#q=WzAsNixbMiwwLCJcXER1YWwgSSJdLFsxLDAsIlxcRHVhbCBKIl0sWzAsMCwiXFxEdWFsXFxjaW1wezF9Il0sWzIsMSwiSV5cXG9wIl0sWzEsMSwiSl5cXG9wIl0sWzAsMSwiXFxjaW1wezF9Xlxcb3AiXSxbMSwwLCJcXER1YWwoZl4qKSJdLFsyLDEsIlxcRHVhbChqXiopIl0sWzIsNSwiPSJdLFsxLDQsIj0iXSxbMCwzLCI9Il0sWzUsNCwial8qXlxcb3AiLDJdLFs0LDMsImZfKl5cXG9wIiwyXV0=
\[\begin{tikzcd}
	{\Dual\cimp{1}} & {\Dual J} & {\Dual I} \\
	{\cimp{1}^\op} & {J^\op} & {I^\op}
	\arrow["{\Dual(j^*)}", from=1-1, to=1-2]
	\arrow["{=}", from=1-1, to=2-1]
	\arrow["{\Dual(f^*)}", from=1-2, to=1-3]
	\arrow["{=}", from=1-2, to=2-2]
	\arrow["{=}", from=1-3, to=2-3]
	\arrow["{j_*^\op}"', from=2-1, to=2-2]
	\arrow["{f_*^\op}"', from=2-2, to=2-3]
\end{tikzcd}\]
in $\Ordi$, where $j^*$ denotes the map $\lad(j_*)\colon \cimp{1}\from
J$ in $\Ordi$, we already have that the outer square commutes as well
as the left one, and hence the desired conclusion.

Let us thus assume that $J = \cimp{1}$.
Then the result follows from the diagram
% https://q.uiver.app/#q=WzAsOCxbMSwxLCJcXER1YWxcXGNpbXB7MX0iXSxbMiwxLCJcXER1YWwgSSJdLFsxLDIsIlxcY2ltcHsxfV5cXG9wIl0sWzIsMiwiSV5cXG9wIl0sWzAsMCwiXFxpZCJdLFszLDAsImZeKiJdLFszLDMsIlxcbWF4KGZeKileey0xfSgwKSJdLFswLDMsIlxcbWF4KFxcY2ltcHsxfV5cXG9wKSA9IDAiXSxbMCwyLCI9Il0sWzEsMywiPSJdLFswLDEsIlxcRHVhbChmXiopIl0sWzIsMywiZl8qXlxcb3AiLDJdLFs0LDAsIlxcaW4iLDMseyJzdHlsZSI6eyJib2R5Ijp7Im5hbWUiOiJub25lIn0sImhlYWQiOnsibmFtZSI6Im5vbmUifX19XSxbNSwxLCJcXGluIiwzLHsic3R5bGUiOnsiYm9keSI6eyJuYW1lIjoibm9uZSJ9LCJoZWFkIjp7Im5hbWUiOiJub25lIn19fV0sWzYsMywiXFxpbiIsMyx7InN0eWxlIjp7ImJvZHkiOnsibmFtZSI6Im5vbmUifSwiaGVhZCI6eyJuYW1lIjoibm9uZSJ9fX1dLFs0LDUsIiIsMyx7InN0eWxlIjp7InRhaWwiOnsibmFtZSI6Im1hcHMgdG8ifX19XSxbNSw2LCIiLDMseyJzdHlsZSI6eyJ0YWlsIjp7Im5hbWUiOiJtYXBzIHRvIn19fV0sWzcsMiwiXFxpbiIsMyx7InN0eWxlIjp7ImJvZHkiOnsibmFtZSI6Im5vbmUifSwiaGVhZCI6eyJuYW1lIjoibm9uZSJ9fX1dLFs0LDcsIiIsMyx7InN0eWxlIjp7InRhaWwiOnsibmFtZSI6Im1hcHMgdG8ifX19XSxbNyw2LCI/IiwwLHsic3R5bGUiOnsidGFpbCI6eyJuYW1lIjoibWFwcyB0byJ9LCJib2R5Ijp7Im5hbWUiOiJkYXNoZWQifX19XV0=
\[\begin{tikzcd}
	\id &&& {f^*} \\
	& {\Dual\cimp{1}} & {\Dual I} \\
	& {\cimp{1}^\op} & {I^\op} \\
	{\max(\cimp{1}^\op) = 0} &&& {\max(f^*)^{-1}(0)}
	\arrow[maps to, from=1-1, to=1-4]
	\arrow["\in"{marking, allow upside down}, draw=none, from=1-1, to=2-2]
	\arrow[maps to, from=1-1, to=4-1]
	\arrow["\in"{marking, allow upside down}, draw=none, from=1-4, to=2-3]
	\arrow[maps to, from=1-4, to=4-4]
	\arrow["{\Dual(f^*)}", from=2-2, to=2-3]
	\arrow["{=}", from=2-2, to=3-2]
	\arrow["{=}", from=2-3, to=3-3]
	\arrow["{f_*^\op}"', from=3-2, to=3-3]
	\arrow["\in"{marking, allow upside down}, draw=none, from=4-1, to=3-2]
	\arrow["{?}", dashed, maps to, from=4-1, to=4-4]
	\arrow["\in"{marking, allow upside down}, draw=none, from=4-4, to=3-3]
\end{tikzcd}\]
in $\Ordt$, where $0 = \min\cimp{1}$, and the computation $f_*(0) =
\max(f^*)^{-1}(0)$ of the value of the right adjoint functor.
\end{proof}

\begin{corollary}
  The diagram
  % https://q.uiver.app/#q=WzAsNixbMCwxLCJcXE9yZF5cXG9wIl0sWzEsMCwiXFxPcmR0Xlxcb3AiXSxbMSwyLCJcXE9yZGleXFxvcCJdLFszLDEsIlxcT3JkIl0sWzIsMCwiXFxPcmRpIl0sWzIsMiwiXFxPcmR0Il0sWzEsNCwiXFxsYWQiXSxbMiw1LCJcXHJhZCIsMl0sWzAsMSwidCJdLFswLDIsImkiLDJdLFs1LDMsIiIsMix7InN0eWxlIjp7InRhaWwiOnsibmFtZSI6Imhvb2siLCJzaWRlIjoidG9wIn19fV0sWzQsMywiIiwyLHsic3R5bGUiOnsidGFpbCI6eyJuYW1lIjoiaG9vayIsInNpZGUiOiJ0b3AifX19XV0=
\[\begin{tikzcd}
	& {\Ordt^\op} & \Ordi \\
	{\Ord^\op} &&& \Ord \\
	& {\Ordi^\op} & \Ordt
	\arrow["\lad", from=1-2, to=1-3]
	\arrow[hook, from=1-3, to=2-4]
	\arrow["t", from=2-1, to=1-2]
	\arrow["i"', from=2-1, to=3-2]
	\arrow["\rad"', from=3-2, to=3-3]
	\arrow[hook, from=3-3, to=2-4]
\end{tikzcd}\]
commutes.
\end{corollary}
\begin{proof}
  This follows from Theorem and Lemma~\ref{lem:functor-to-pull-back}.
\end{proof}

\begin{notation}
  $[-]\colon \Ord^\op\to \Ordit = \Ordi\times_\Ord \Ordt$
  (Lemma~\ref{lem:pull-back}) denotes the unique functor to make the
  diagram
  % https://q.uiver.app/#q=WzAsNyxbMCwxLCJcXE9yZF5cXG9wIl0sWzEsMCwiXFxPcmR0Xlxcb3AiXSxbMSwyLCJcXE9yZGleXFxvcCJdLFszLDEsIlxcT3JkIl0sWzIsMCwiXFxPcmRpIl0sWzIsMiwiXFxPcmR0Il0sWzEsMSwiXFxPcmRpdCJdLFsxLDQsIlxcbGFkIl0sWzIsNSwiXFxyYWQiLDJdLFswLDEsInQiXSxbMCwyLCJpIiwyXSxbNSwzLCIiLDIseyJzdHlsZSI6eyJ0YWlsIjp7Im5hbWUiOiJob29rIiwic2lkZSI6InRvcCJ9fX1dLFs0LDMsIiIsMix7InN0eWxlIjp7InRhaWwiOnsibmFtZSI6Imhvb2siLCJzaWRlIjoidG9wIn19fV0sWzYsNCwiIiwyLHsic3R5bGUiOnsidGFpbCI6eyJuYW1lIjoiaG9vayIsInNpZGUiOiJ0b3AifX19XSxbNiw1LCIiLDIseyJzdHlsZSI6eyJ0YWlsIjp7Im5hbWUiOiJob29rIiwic2lkZSI6InRvcCJ9fX1dLFswLDYsIlstXSJdXQ==
\[\begin{tikzcd}
	& {\Ordt^\op} & \Ordi \\
	{\Ord^\op} & \Ordit && \Ord \\
	& {\Ordi^\op} & \Ordt
	\arrow["\lad", from=1-2, to=1-3]
	\arrow[hook, from=1-3, to=2-4]
	\arrow["t", from=2-1, to=1-2]
	\arrow["{[-]}", from=2-1, to=2-2]
	\arrow["i"', from=2-1, to=3-2]
	\arrow[hook, from=2-2, to=1-3]
	\arrow[hook, from=2-2, to=3-3]
	\arrow["\rad"', from=3-2, to=3-3]
	\arrow[hook, from=3-3, to=2-4]
\end{tikzcd}\]
commute (see Corollary).
For an object $I$ of $\Ord$, $[-](I)$ is denoted by $[I]$.
\end{notation}

\begin{lemma}
  The diagram
  % https://q.uiver.app/#q=WzAsNCxbMCwwLCJcXE9yZF5cXG9wIl0sWzAsMSwiXFxPcmReXFxvcCJdLFsxLDAsIlxcT3JkaXQiXSxbMSwxLCJcXE9yZGl0Il0sWzAsMiwiWy1dIl0sWzEsMywiWy1dIiwyXSxbMCwxLCJcXG9wIiwyXSxbMiwzLCJcXG9wIl0sWzAsMywiXFxCaXIiXV0=
\[\begin{tikzcd}
	{\Ord^\op} & \Ordit \\
	{\Ord^\op} & \Ordit
	\arrow["{[-]}", from=1-1, to=1-2]
	\arrow["\op"', from=1-1, to=2-1]
	\arrow["\Bir", from=1-1, to=2-2]
	\arrow["\op", from=1-2, to=2-2]
	\arrow["{[-]}"', from=2-1, to=2-2]
\end{tikzcd}\]
commutes.
\end{lemma}
\begin{proof}
  In view of the description of the Birkhoff functor as a functor to a
  fibre product using Lemma~\ref{lem:functor-to-pull-back}, the
  commutativity of the upper triangle follows from the diagram
  % https://q.uiver.app/#q=WzAsMTAsWzAsMSwiXFxPcmReXFxvcCJdLFsxLDIsIlxcT3JkaV5cXG9wIl0sWzEsMCwiXFxPcmR0Xlxcb3AiXSxbMiwyLCJcXE9yZHQiXSxbMiwwLCJcXE9yZGkiXSxbMSwxLCJcXE9yZGl0Il0sWzMsMiwiXFxPcmRpIl0sWzQsMSwiXFxPcmQiXSxbMywwLCJcXE9yZHQiXSxbMiwxLCJcXE9yZGl0Il0sWzEsMywiXFxyYWQiXSxbMiw0LCJcXGxhZCIsMl0sWzAsMSwiaSIsMl0sWzAsMiwidCJdLFs1LDMsIiIsMix7InN0eWxlIjp7InRhaWwiOnsibmFtZSI6Imhvb2siLCJzaWRlIjoidG9wIn19fV0sWzUsNCwiIiwyLHsic3R5bGUiOnsidGFpbCI6eyJuYW1lIjoiaG9vayIsInNpZGUiOiJ0b3AifX19XSxbMCw1LCJbLV0iXSxbMyw2LCJcXG9wIl0sWzYsNywiIiwwLHsic3R5bGUiOnsidGFpbCI6eyJuYW1lIjoiaG9vayIsInNpZGUiOiJ0b3AifX19XSxbOCw3LCIiLDAseyJzdHlsZSI6eyJ0YWlsIjp7Im5hbWUiOiJob29rIiwic2lkZSI6InRvcCJ9fX1dLFs0LDgsIlxcb3AiLDJdLFs1LDksIlxcb3AiXSxbOSw2LCIiLDAseyJzdHlsZSI6eyJ0YWlsIjp7Im5hbWUiOiJob29rIiwic2lkZSI6InRvcCJ9fX1dLFs5LDgsIiIsMCx7InN0eWxlIjp7InRhaWwiOnsibmFtZSI6Imhvb2siLCJzaWRlIjoidG9wIn19fV0sWzEsNiwiXFxEdWFsIiwyLHsiY3VydmUiOjN9XSxbMiw4LCJcXER1YWwiLDAseyJjdXJ2ZSI6LTN9XV0=
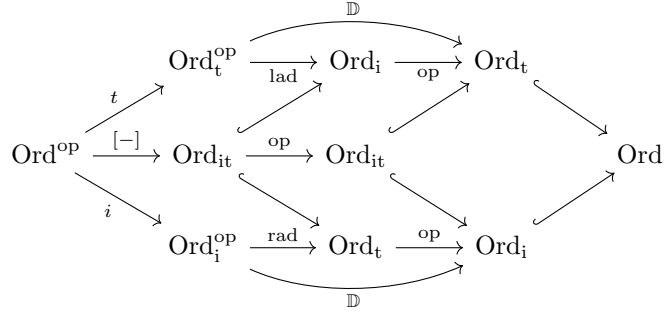
\[\begin{tikzcd}
	& {\Ordt^\op} & \Ordi & \Ordt \\
	{\Ord^\op} & \Ordit & \Ordit && \Ord \\
	& {\Ordi^\op} & \Ordt & \Ordi
	\arrow["\lad"', from=1-2, to=1-3]
	\arrow["\Dual", curve={height=-18pt}, from=1-2, to=1-4]
	\arrow["\op"', from=1-3, to=1-4]
	\arrow[hook, from=1-4, to=2-5]
	\arrow["t", from=2-1, to=1-2]
	\arrow["{[-]}", from=2-1, to=2-2]
	\arrow["i"', from=2-1, to=3-2]
	\arrow[hook, from=2-2, to=1-3]
	\arrow["\op", from=2-2, to=2-3]
	\arrow[hook, from=2-2, to=3-3]
	\arrow[hook, from=2-3, to=1-4]
	\arrow[hook, from=2-3, to=3-4]
	\arrow["\rad", from=3-2, to=3-3]
	\arrow["\Dual"', curve={height=18pt}, from=3-2, to=3-4]
	\arrow["\op", from=3-3, to=3-4]
	\arrow[hook, from=3-4, to=2-5]
\end{tikzcd}\]
which commutes by Theorem~\ref{thm:adjoint-dual}.

The commutativity of the lower triangle follows from this and the
commutation of $\Bir$ and $\op$, or from a similar diagram as above.
\end{proof}

\begin{remark}
  This could be the definition of $[-]$.
\end{remark}

From Lemma, Theorem~\ref{thm:adjoint-dual} and
Proposition~\ref{prop:sequence-duality}, we obtain the following
result.

\begin{proposition}
  The diagram
  % https://q.uiver.app/#q=WzAsOCxbMCwxLCJcXE9yZF5cXG9wIl0sWzIsMCwiXFxPcmR0Xlxcb3AiXSxbMiwyLCJcXE9yZGleXFxvcCJdLFs0LDEsIlxcT3JkaXReXFxvcCJdLFsyLDEsIlxcT3JkaXQiXSxbNiwxLCJcXE9yZCJdLFs0LDAsIlxcT3JkaSJdLFs0LDIsIlxcT3JkdCJdLFswLDQsIlstXSJdLFszLDUsIlstXV57LTF9IiwyXSxbMSw2LCJcXGxhZCJdLFsyLDcsIlxccmFkIiwyXSxbMCwxLCJ0Il0sWzEsMywiaSIsMix7ImxhYmVsX3Bvc2l0aW9uIjo3MH1dLFswLDIsImkiLDJdLFsyLDMsInQiLDAseyJsYWJlbF9wb3NpdGlvbiI6NzB9XSxbNCw2LCIiLDIseyJzdHlsZSI6eyJ0YWlsIjp7Im5hbWUiOiJob29rIiwic2lkZSI6InRvcCJ9fX1dLFs0LDcsIiIsMix7InN0eWxlIjp7InRhaWwiOnsibmFtZSI6Imhvb2siLCJzaWRlIjoidG9wIn19fV0sWzcsNSwiIiwyLHsic3R5bGUiOnsidGFpbCI6eyJuYW1lIjoiaG9vayIsInNpZGUiOiJ0b3AifX19XSxbNiw1LCIiLDIseyJzdHlsZSI6eyJ0YWlsIjp7Im5hbWUiOiJob29rIiwic2lkZSI6InRvcCJ9fX1dXQ==
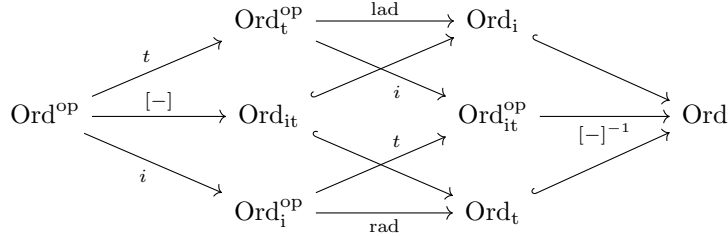
\[\begin{tikzcd}
	&& {\Ordt^\op} && \Ordi \\
	{\Ord^\op} && \Ordit && {\Ordit^\op} && \Ord \\
	&& {\Ordi^\op} && \Ordt
	\arrow["\lad", from=1-3, to=1-5]
	\arrow["i"'{pos=0.7}, from=1-3, to=2-5]
	\arrow[hook, from=1-5, to=2-7]
	\arrow["t", from=2-1, to=1-3]
	\arrow["{[-]}", from=2-1, to=2-3]
	\arrow["i"', from=2-1, to=3-3]
	\arrow[hook, from=2-3, to=1-5]
	\arrow[hook, from=2-3, to=3-5]
	\arrow["{[-]^{-1}}"', from=2-5, to=2-7]
	\arrow["t"{pos=0.7}, from=3-3, to=2-5]
	\arrow["\rad"', from=3-3, to=3-5]
	\arrow[hook, from=3-5, to=2-7]
\end{tikzcd}\]
commutes, giving an isomorphism of pull-back squares.
\end{proposition}

\section{Algebraic operations}
\label{sec:operations}

In this section, we study how the Birkhoff-type dualities interact
with some algebraic structures on the relevant categories.

\subsection{Universal modules over the universal monoid}
\label{sec:ord-module}
In this section, we study some simple algebraic structures on some of
our categories.

In order to fix our conventions, let us begin with explicitly giving
a (perhaps not so common) formulation of the definition of some
familiar notions.

In our terminology, we let \kore{associativity} to include unitality
unless otherwise stated.

\begin{definition}
  \label{def:enriched-category}
  Let $\cat{V}$ be a monoidal category with monoidal operation
  $\tensor$.
  Then a \kore{category} $\cat{C}$ \kore{enriched} in $\cat{V}$
  consists of the following pieces of information.
  \begin{itemize}
  \item A collection $\Ob\cat{C}$ of ``\kore{objects}'' of $\cat{C}$.
  \item For every pair $X$ of objects of $\cat{C}$, an object
    $\Hom_{\cat{C}}[X]$ of $\cat{V}$.
  \item Suppose given the following data.
    \begin{itemize}
    \item An object $I$ of $\Ord$.
    \item A family $X = (X_b)_{b\in \Bir I}$ of objects of $\cat{C}$
      indexed by $\Bir I$.
    \end{itemize}
    Then a ``\kore{composition}'' morphism
    \[
      \Compose_I\colon
      \Hom_{\cat{C}}[X]\longto\Hom_{\cat{C}}[\boundary X],
    \]
    where 
    \begin{itemize}
    \item $\Hom_{\cat{C}}[X]
      := \Tensor_{i\in I(= \Bir \Bir I)}\Hom_{\cat{C}}[X\resto{i}]$,
      where $X\resto{i}$ for the morphism $i\colon \Bir I\to \cimp{1}$
      in $\Ordit$ is the pair $(X_{\max i^{-1}(\min\cimp{1})},
      X_{\min i^{-1}(\max\cimp{1})})$, and
    \item $\boundary X$ is the pair $(X_{\min \Bir I}, X_{\max \Bir I})$.
    \end{itemize}
  \end{itemize}
  Moreover, the composition is required to be associative.
\end{definition}

We write $\Hom_{\cat{C}}(X, Y)$ for $\Hom_{\cat{C}}[(X, Y)]$ (and not
$\Hom_{\cat{C}}[(Y, X)]$, but this choice is inessential).
The use of $\Bir I$ rather than $[I]$ was only for convenience, and
this choice is not essential either.

A monoid in $\cat{V}$ is essentially a special case of this.

\begin{definition}
  Given an object $A$ of $\cat{V}$, the structure of a \kore{monoid}
  on $A$ consists of the association to every object
  $I$ of $\Ord$ of a ``\kore{multiplication}'' morphism
  \[
    \Tensor_I\colon A^{\tensor I}\longto A
  \]
  such that the following data defines a category $BA$ enriched
  in $\cat{V}$.
  \begin{itemize}
  \item Any non-empty collection $O$ as $\Ob BA$.
  \item $A$ as $\Hom_{BA}[X]$ for every pair $X$ in $O$.
  \item Suppose given the following data.
    \begin{itemize}
    \item An object $I$ of $\Ord$.
    \item A family $X$ indexed by $\Bir I$ in $O$.
    \end{itemize}
    Then $\Tensor_I$ as the composition functor $\Hom_{BA}[X]
    \to \Hom_{BA}[\boundary X]$.
  \end{itemize}
\end{definition}

\begin{example}
  Let $\cat{V}$ be a monoidal category, $\cat{C}$ a category enriched
  in $\cat{V}$, and $X$ an object of $\cat{C}$.
  Then the object $\End_{\cat{C}}(X) = \Hom_{\cat{C}}(X, X)$ of
  $\cat{V}$ has the structure of a monoid whose multiplication is
  given by the composition in $\cat{C}$.
\end{example}

\begin{example}
  In the case where $\cat{V}$ is the Cartesian closed category
  $\Catset$, the structure of a monoid on a object of $X$ of $\cat{V}$
  is equivalent to a monoidal structure on the category $X$.
\end{example}

Let us give $\Ord$ a monoidal structure.
Given an object $I$ of $\Ord$, we would like to define the
multiplication functor
\[
  \Disj_I\colon \Ord^{\times I}\longto \Ord.
\]
For a family $X = (X_i)_{i\in I}$ of objects of $\Ord$, we let
$\Disj_I X$ be the object which accepts a morphism to $I$ which
is characterized by the condition that, for every $i\in I$, we
have a pull-back square
% https://q.uiver.app/#q=WzAsNCxbMSwwLCJcXERpc2pfSSBYIl0sWzAsMSwiKiJdLFswLDAsIlhfaSJdLFsxLDEsIkkiXSxbMiwwXSxbMiwxXSxbMSwzLCJpIl0sWzAsM11d
\[\begin{tikzcd}
	{X_i} & {\Disj_I X} \\
	{*} & I
	\arrow[from=1-1, to=1-2]
	\arrow[from=1-1, to=2-1]
	\arrow[from=1-2, to=2-2]
	\arrow["i", from=2-1, to=2-2]
\end{tikzcd}\]
where $*$ denotes the terminal set.

If $\cat{V}$ is a \emph{symmetric} monoidal category, then we can
define the opposite of a $\cat{V}$-enriched category, and in
particular, of a monoid in $\cat{V}$.
Our monoidal structure makes $\Ord$ a self-opposite monoid.
Indeed, the functor $\op\colon \Ord\to \Ord$ is an isomorphism of
the monoid $\Ord$ to its monoidal opposite.

\begin{remark}
  \label{rem:universal-monoid}
  The described monoidal structure of $\Ord$ makes its terminal object
  a universal monoid, so monoids in a monoidal category $\cat{V}$ are
  classifed by monoidal functors $\Ord\to \cat{V}$.
\end{remark}

\begin{notation}
  $\Univ$ denotes the terminal monoid in $\Ord$.
\end{notation}

\begin{terminology}
  \label{term:module}
  Let $\cat{V}$ be a monoidal category, and let $A$ be a monoid in
  $\cat{V}$.
  Then for a category $\cat{C}$ enriched in $\cat{V}$, a
  \kore{module over $A$}, or an \kore{$A$-module} in $\cat{C}$ is a
  functor $BA\to \cat{C}$ of $\cat{V}$-enriched categories.
\end{terminology}

Equivalently, an $A$-module in $\cat{C}$ is an object $M$ of $\cat{C}$
equipped with a morphism $A\to \End_{\cat{C}}(M)$ of monoids in
$\cat{V}$.

\begin{terminology}
  Let $\cat{V}$ be a \emph{symmetric} monoidal category, and let $A$
  be an associative monoid in $\cat{V}$.
  Then for a category $\cat{C}$ enriched in $\cat{V}$, a
  \kore{contravariant module} over $A$ in $\cat{C}$ is a module in
  $\cat{C}$ over the opposite of the monoid $A$.
\end{terminology}

Note that a contravariant $A$-module, despite its name, is \emph{not}
a module over $A$ in the previous sense.
For distinction, we introduce another term.

\begin{terminology}
  Let $\cat{V}$ be a symmetric monoidal category, and let $A$ be an
  associative monoid in $\cat{V}$.
  Then for a category $\cat{C}$ enriched in $\cat{V}$, a
  \kore{covariant} module over $A$ in $\cat{C}$ is a module over
  $A$ in $\cat{C}$ in the sense of Terminology~\ref{term:module}.
\end{terminology}

When $\cat{V}$ is symmetric, the term ``\kore{module}'' may more often
mean either a covariant or a contravariant module than only a
covariant one.

If $\cat{V}$ is symmetric closed, then it is enriched in itself, so,
we can speak of both a covariant module and a contravariant module in
$\cat{V}$.

The monoidal structure of $\Ord$ of course makes it a bimodule over
itself (with commuting co- and contra- variant actions, assuming that
we do not need to call it a ``bivariant'' module).
$\op\colon \Ord\to \Ord$ is \emph{not} a morphism of bimodules but
is compatible with the isomorphism $\op$ on the scalars.

The subcategory $\Ordt$ of $\Ord$ is closed under the \emph{co}variant
action on $\Ord$ by $\Ord$.
It follows that $\Ordt$ has the structure of a covariant module over
$\Ord$ determined by the requirement that the forgetful functor
$\Ordt\into \Ord$ be a homomorphism of covariant $\Ord$-modules.

\begin{remark}
  \label{rem:universal-covariant-module}
  Given a monoid $A$ in a monoidal category $\cat{V}$, there is a
  notion of \kore{$A$-module in a covariant $\cat{V}$-module}
  $\cat{M}$, where $\cat{M}$ is in fact such an object in a
  categorified sense,
  so some, though not us, might prefer calling it something like ``a
  $2$-module'' to emphasize the enrichment.
  (Unfortunately, the term ``bimodule'' has already been given a
  different meaning.)
  An $A$-module in $\cat{M}$ amounts to a lax morphism from the
  terminal covariant $\Ord$-module to $A^*\cat{M}$, where the target
  is the covariant $\Ord$-module obtained from $\cat{M}$ through the
  monoidal functor $\Ord\to \cat{V}$ classifying $A$.
  Then the terminal (in fact, zero) object in the covariant
  $\Ord$-module $\Ordt$ is universal among modules over the universal
  monoid $\Univ$ in a covariant $\Ord$-module.
  It follows that, given a monoid $A$ in a monoidal category
  $\cat{V}$, $A$-modules in a covariant $\cat{V}$-module $\cat{M}$ are
  classified by morphisms $\Ordt\to A^*\cat{M}$ of covariant
  $\Ord$-modules.
  (In terms of a monad, $A^*\cat{M}$ may be described as the
  underlying category of $\cat{M}$ equipped with the monad to be the
  free $A$-module monad.
  See Remark~\ref{rem:universal-monad-algebra}.)
\end{remark}

The forgetful functor $\Ordi\into \Ord$ is similarly the inclusion of
a \emph{contra}variant submodule of $\Ord$.
This submodule is exchanged with the covariant submodule $\Ordt$ by
the functor $\op\colon \Ord\to \Ord$, which flips the variance
of the modules.

\begin{remark}
  \label{rem:universal-contravariant-module}
  The zero object in the contravariant $\Ord$-module $\Ordi$ is
  universal among modules over the monoid $\Univ$ in a contravariant
  $\Ord$-module.
\end{remark}

\begin{theorem}\label{thm:adjoint-homomorphism}
  The functor $i\colon \Ord\to \Ordi$ lifts to a left adjoint of the
  forgetful functor $\Ordi\into \Ord$ in the bicategory of
  contravariant $\Ord$-modules.
\end{theorem}
\begin{proof}
  Given an object $A$ of $\Ord$, let the contravariant action of $A$
  on the underlying category of an $\Ord$-module as a functor be
  denoted by $- \disj A$.
  Then, for every object $X$ of $\Ord$, we have a map $a\colon
  i(X \disj A)\to iX \disj A$ in $\Ordi$ which makes the diagram
  % https://q.uiver.app/#q=WzAsMyxbMCwxLCJpKFhcXGRpc2ogQSkiXSxbMSwwLCJYXFxkaXNqIEEiXSxbMiwxLCJpWFxcZGlzaiBBIl0sWzAsMiwiYSIsMl0sWzEsMiwiXFxldGFcXGRpc2ogQSJdLFsxLDAsIlxcZXRhIiwyXV0=
\[\begin{tikzcd}
	& {X\disj A} \\
	{i(X\disj A)} && {iX\disj A}
	\arrow["\eta"', from=1-2, to=2-1]
	\arrow["{\eta\disj A}", from=1-2, to=2-3]
	\arrow["a"', from=2-1, to=2-3]
\end{tikzcd}\]
in $\Ord$ commute, namely, the map corresponding to $\eta\disj A$
under the adjunction in $\Cat$ between $i$ and the forgetful functor.
It suffices to show that $a$ is an isomorphism for every $X$ and $A$,
but this is easy to see by inspection.
\end{proof}

\begin{corollary}
  \label{cor::adjoint-homomorphism}
  The functor $t\colon \Ord\to \Ordt$ lifts to a left adjoint of the
  forgetful functor $\Ordt\into \Ord$ in the bicategory of
  covariant $\Ord$-modules.
\end{corollary}
\begin{proof}
  Apply $\op$ to Theorem, or apply the similar arguments.
\end{proof}

\subsection{Pairing}
\label{sec:pairing}

In this section, we introduce some further algebraic structure, which
is a pairing between the modules we have considered above.

\begin{notation}
  Given an object $I$ of $\Ord$ and a family $X = (X_i)_{i\in I}$ of
objects of a monoidal category $\cat{V}$, $X_{\min I}\tensor
\cdots\tensor X_{\max I}$, where the factors are written in the order
of the indices in $I$, denotes the monoidal product $\Tensor_I X$.
\end{notation}

Just in order to motivate the conventions chosen in
Notation~\ref{nota:module-action} below, let us discuss some
heuristics.
One can extend the definition of a \kore{co-} (resp.~\kore{contra-})
\kore{variant module} over a monoid $A$ in $\cat{V}$ to the case where
$\cat{V}$ is not closed but merely associative.
It is just a module over $A$ in the \emph{co-} (resp.~\emph{contra-})
variant $\cat{V}$-module $\cat{V}$, as considered in
Remark~\ref{rem:universal-covariant-module}
(resp.~\ref{rem:universal-contravariant-module}).
This will be a left (resp.~right) $A$-module in $\cat{V}$ with
respect to the arbitrarily chosen notation above.

\begin{remark}
  The definitions make sense more generally in any bimodule $\cat{M}$
  over $\cat{V}$.
  However, if the monoidal structure of $\cat{V}$ is given a symmetry,
  then one needs some restriction on the bimodule structure to have
  that a contravariant $A$-module be the same as a covariant module
  over the opposite of the monoid $A$.
  It suffices if $\cat{M}$ comes from a module over the symmetric
  monoidal category.
  The question would be subtler when $\cat{V}$ is given only a
  braiding, but it will suffice again if $\cat{M}$ comes from an
  operadic $\cat{V}$-module for the little \emph{disks} operad.
  In particular, the bimodule $\cat{V}$ is fine in either case.
\end{remark}

\begin{notation}
  \label{nota:module-action}
  $\Disj\colon \Ord\times\Ordt\to \Ordt$ (resp.~$\Disj\colon
  \Ordi\times\Ord\to \Ordi$) denotes the structure action of $\Ord$ on
  its co- (resp.~contra-) variant module.
\end{notation}

With these choices of sides, we would hopefully not need to use the
symmetry of the (Cartesian) monoidal structure on $\Catset$ more often
than otherwise.

\begin{proposition}
  The diagram
  \begin{equation*}
    % https://q.uiver.app/#q=WzAsNixbMCwxLCJcXE9yZGlcXHRpbWVzIFxcT3JkdCJdLFsxLDIsIlxcT3JkaVxcdGltZXNcXE9yZCJdLFsyLDIsIlxcT3JkaSJdLFszLDEsIlxcT3JkIl0sWzEsMCwiXFxPcmRcXHRpbWVzIFxcT3JkdCJdLFsyLDAsIlxcT3JkdCJdLFswLDEsIiIsMix7InN0eWxlIjp7InRhaWwiOnsibmFtZSI6Imhvb2siLCJzaWRlIjoidG9wIn19fV0sWzAsNCwiIiwwLHsic3R5bGUiOnsidGFpbCI6eyJuYW1lIjoiaG9vayIsInNpZGUiOiJ0b3AifX19XSxbNCw1LCJcXERpc2oiXSxbMSwyLCJcXERpc2oiLDJdLFsyLDMsIiIsMix7InN0eWxlIjp7InRhaWwiOnsibmFtZSI6Imhvb2siLCJzaWRlIjoidG9wIn19fV0sWzUsMywiIiwwLHsic3R5bGUiOnsidGFpbCI6eyJuYW1lIjoiaG9vayIsInNpZGUiOiJ0b3AifX19XV0=
\begin{tikzcd}
	& {\Ord\times \Ordt} & \Ordt \\
	{\Ordi\times \Ordt} &&& \Ord \\
	& {\Ordi\times\Ord} & \Ordi
	\arrow["\Disj", from=1-2, to=1-3]
	\arrow[hook, from=1-3, to=2-4]
	\arrow[hook, from=2-1, to=1-2]
	\arrow[hook, from=2-1, to=3-2]
	\arrow["\Disj"', from=3-2, to=3-3]
	\arrow[hook, from=3-3, to=2-4]
\end{tikzcd}
  \end{equation*}
  commutes, where each of the unlabeled arrows is (induced by) the
  relevant fogetful functor.
\end{proposition}
\begin{proof}
  This follows from the diagram
  % https://q.uiver.app/#q=WzAsNyxbMCwxLCJcXE9yZGlcXHRpbWVzIFxcT3JkdCJdLFsxLDIsIlxcT3JkaVxcdGltZXNcXE9yZCJdLFsyLDIsIlxcT3JkaSJdLFszLDEsIlxcT3JkIl0sWzEsMCwiXFxPcmRcXHRpbWVzIFxcT3JkdCJdLFsyLDAsIlxcT3JkdCJdLFsyLDEsIlxcT3JkXFx0aW1lcyBcXE9yZCJdLFswLDEsIiIsMix7InN0eWxlIjp7InRhaWwiOnsibmFtZSI6Imhvb2siLCJzaWRlIjoidG9wIn19fV0sWzAsNCwiIiwwLHsic3R5bGUiOnsidGFpbCI6eyJuYW1lIjoiaG9vayIsInNpZGUiOiJ0b3AifX19XSxbNCw1LCJcXERpc2oiXSxbMSwyLCJcXERpc2oiLDJdLFsyLDMsIiIsMix7InN0eWxlIjp7InRhaWwiOnsibmFtZSI6Imhvb2siLCJzaWRlIjoidG9wIn19fV0sWzUsMywiIiwwLHsic3R5bGUiOnsidGFpbCI6eyJuYW1lIjoiaG9vayIsInNpZGUiOiJ0b3AifX19XSxbNiwzLCJcXERpc2oiXSxbMSw2LCIiLDAseyJzdHlsZSI6eyJ0YWlsIjp7Im5hbWUiOiJob29rIiwic2lkZSI6InRvcCJ9fX1dLFs0LDYsIiIsMix7InN0eWxlIjp7InRhaWwiOnsibmFtZSI6Imhvb2siLCJzaWRlIjoidG9wIn19fV1d
\[\begin{tikzcd}
	& {\Ord\times \Ordt} & \Ordt \\
	{\Ordi\times \Ordt} && {\Ord\times \Ord} & \Ord \\
	& {\Ordi\times\Ord} & \Ordi
	\arrow["\Disj", from=1-2, to=1-3]
	\arrow[hook, from=1-2, to=2-3]
	\arrow[hook, from=1-3, to=2-4]
	\arrow[hook, from=2-1, to=1-2]
	\arrow[hook, from=2-1, to=3-2]
	\arrow["\Disj", from=2-3, to=2-4]
	\arrow[hook, from=3-2, to=2-3]
	\arrow["\Disj"', from=3-2, to=3-3]
	\arrow[hook, from=3-3, to=2-4]
\end{tikzcd}\]
  which commutes since the forgetful functor
  $\Ordt \text{(resp.~$\Ordi$)}\to \Ord$ is a homomorphism of co-
  (resp.~contra-) variant $\Ord$-modules.
\end{proof}

\begin{notation}
  $\Disj\colon \Ordi\times \Ordt\to \Ordit$ denotes the functor
  determined from the commutative diagram of Proposition by the
  pull-back square Lemma~\ref{lem:pull-back}.
\end{notation}

Thus, we have obtained an $\Ordit$-valued pairing on between $\Ordi$
and $\Ordt$.

\begin{theorem}
  The functor $\Disj\colon \Ordi\times \Ordt\to \Ordit$ coequalizes
  the two functors
  % https://q.uiver.app/#q=WzAsMixbMCwwLCJcXE9yZGlcXHRpbWVzXFxPcmRcXHRpbWVzXFxPcmR0Il0sWzEsMCwiXFxPcmRpXFx0aW1lc1xcT3JkdCJdLFswLDEsIlxcaWRcXHRpbWVzIFxcRGlzaiIsMCx7Im9mZnNldCI6LTF9XSxbMCwxLCJcXERpc2pcXHRpbWVzIFxcaWQiLDIseyJvZmZzZXQiOjF9XV0=
\[\begin{tikzcd}
	{\Ordi\times\Ord\times\Ordt} & {\Ordi\times\Ordt}
	\arrow["{\id\times \Disj}", shift left, from=1-1, to=1-2]
	\arrow["{\Disj\times \id}"', shift right, from=1-1, to=1-2]
\end{tikzcd}\]
\end{theorem}
\begin{proof}
  The commutative diagram
  % https://q.uiver.app/#q=WzAsNyxbMCwxLCJcXE9yZGlcXHRpbWVzXFxPcmRcXHRpbWVzXFxPcmR0Il0sWzEsMSwiXFxPcmRpXFx0aW1lc1xcT3JkXFx0aW1lc1xcT3JkIl0sWzIsMCwiXFxPcmRpXFx0aW1lc1xcT3JkIl0sWzMsMSwiXFxPcmRpIl0sWzEsMCwiXFxPcmRpXFx0aW1lc1xcT3JkdCJdLFsyLDIsIlxcT3JkaVxcdGltZXNcXE9yZCJdLFsxLDIsIlxcT3JkaVxcdGltZXNcXE9yZHQiXSxbMiwzLCJcXERpc2oiXSxbMCwxLCIiLDAseyJzdHlsZSI6eyJ0YWlsIjp7Im5hbWUiOiJob29rIiwic2lkZSI6InRvcCJ9fX1dLFsxLDIsIlxcaWRcXHRpbWVzIFxcRGlzaiIsMl0sWzAsNCwiXFxpZFxcdGltZXMgXFxEaXNqIl0sWzQsMiwiIiwwLHsic3R5bGUiOnsidGFpbCI6eyJuYW1lIjoiaG9vayIsInNpZGUiOiJ0b3AifX19XSxbMSw1LCJcXERpc2pcXHRpbWVzIFxcaWQiXSxbNSwzLCJcXERpc2oiLDJdLFswLDYsIlxcRGlzalxcdGltZXMgXFxpZCIsMl0sWzYsNSwiIiwyLHsic3R5bGUiOnsidGFpbCI6eyJuYW1lIjoiaG9vayIsInNpZGUiOiJ0b3AifX19XV0=
  \[\begin{tikzcd}
      & {\Ordi\times\Ordt} & {\Ordi\times\Ord} \\
      {\Ordi\times\Ord\times\Ordt} & {\Ordi\times\Ord\times\Ord} && \Ordi \\
      & {\Ordi\times\Ordt} & {\Ordi\times\Ord}
      \arrow[hook, from=1-2, to=1-3]
      \arrow["\Disj", from=1-3, to=2-4]
      \arrow["{\id\times \Disj}", from=2-1, to=1-2]
      \arrow[hook, from=2-1, to=2-2]
      \arrow["{\Disj\times \id}"', from=2-1, to=3-2]
      \arrow["{\id\times \Disj}"', from=2-2, to=1-3]
      \arrow["{\Disj\times \id}", from=2-2, to=3-3]
      \arrow[hook, from=3-2, to=3-3]
      \arrow["\Disj"', from=3-3, to=2-4]
    \end{tikzcd}\]
  show that the two functors in question are coequalized by the
  restriction of $\Disj\colon \Ordi\times \Ord\to \Ordi$ to
  $\Ordi\times \Ordt$.
  By symmetry or the similar diagram, the two functors are also
  coequalized by the restriction of $\Disj\colon \Ord\times \Ordt\to
  \Ordt$ to $\Ordi\times \Ordt$.
  The result follows from these and Lemma~\ref{lem:pull-back}.

  Alternatively, the restriction of $\Disj\colon \Ordi\times \Ord\to
  \Ordi$ to $\Ordi\times \Ordt$ factorizes as
  % https://q.uiver.app/#q=WzAsMyxbMSwwLCJcXE9yZGl0Il0sWzAsMCwiXFxPcmRpXFx0aW1lc1xcT3JkdCJdLFsyLDAsIlxcT3JkaSJdLFsxLDAsIlxcRGlzaiJdLFswLDIsIiIsMCx7InN0eWxlIjp7InRhaWwiOnsibmFtZSI6Imhvb2siLCJzaWRlIjoidG9wIn19fV1d
  \[\begin{tikzcd}
      {\Ordi\times\Ordt} & \Ordit & \Ordi
      \arrow["\Disj", from=1-1, to=1-2]
      \arrow[hook, from=1-2, to=1-3]
    \end{tikzcd}\]
  and the forgetful functor in this composition is a monomorphism.
\end{proof}

\begin{corollary}
  \label{cor:pairing-adjoint-yuru}
  The functor $\Ordt\to \Hom_{\Cat}(\Ordi, \Ordit)$ corresponding to
  the pairing $\Disj\colon \Ordi\times \Ordt\to \Ordit$ is a
  homomorphism of covariant $\Ord$-modules.
  Equivalently, the functor $\Ordi\to \Hom_{\Cat}(\Ordt, \Ordit)$
  corresponding to the pairing $\Disj$ is a homomorphism of
  contravariant $\Ord$-modules.
\end{corollary}

\subsection{Interactions of the dual structures}
\label{sec:dual-structure}

In this section, we give some concluding results.

The monoidal structure on $\Ord$ can of course be carried over along
an isomorphism $\Ord^\op\equivwith \Ordit$ to make the isomorphism
that of monoids.
We obtain two monoidal structures on $\Ordit$ by using the
isomorphisms $\Bir$ and $[-]$.
These monoidal structures are opposite to each other, and the
resulting two monoids are isomorphic by the automorphism $\op$ of the
category $\Ordit$.

Explicitly, the monoidal structure on $\Ordit$ obtained through $\Bir$
(resp.~$[-]$) can be described as follows.
Given an object $I$ of $\Ord$, the multiplication functor
\[
  \sum_I\colon \Ordit^{\times I}\longto \Ordit
\]
is as follows.
For a family $X = (X_i)_{i\in I}$ of objects of $\Ordit$, $\sum_I X$
is the object which accepts a morphism from $\Bir I$ (resp.~$[I]$)
which is characterized by the condition that, for every $i\in I$, we
have a push-out square
% https://q.uiver.app/#q=WzAsNCxbMCwwLCJcXHRleHR7JFxcQmlyIEkkIChyZXNwLn4kW0ldJCl9Il0sWzEsMCwiXFxzdW1fSSBYIl0sWzAsMSwiXFxjaW1wezF9Il0sWzEsMSwiWF9pIl0sWzAsMiwiaSIsMl0sWzAsMV0sWzEsM10sWzIsM11d
\[\begin{tikzcd}
	{\text{$\Bir I$ (resp.~$[I]$)}} & {\sum_I X} \\
	{\cimp{1}} & {X_i}
	\arrow[from=1-1, to=1-2]
	\arrow["i"', from=1-1, to=2-1]
	\arrow[from=1-2, to=2-2]
	\arrow[from=2-1, to=2-2]
\end{tikzcd}\]
where $i$ is interpreted as a morphism
$\Bir I \text{(resp.~$[I] = \Bir I^\op$)}\to \cimp{1}$ in $\Ordit$.

\kore{We let $\sum$ be the monoidal structure on $\Ordit$ obtained
  through the isomorphism $\Bir$ from $\Disj$ on $\Ord^\op$.}

\begin{remark}
  A meaning of this monoidal structure to universal algebra is that
  the initial comonoid in $\Ordit$ with respect to this structure,
  which is $\Bir\Univ$, is universal among all comonoids in a monoidal
  category.
\end{remark}

Using the isomorphism $\Bir\colon \Ordit\to \Ord^\op$ of
\emph{monoids}, we can now consider $\Ord^\op$ as a bimodule over
$\Ordit$.
This bimodule has a co- (resp.~contra-) variant submodule $\Ordt^\op$
(resp.~$\Ordi^\op$).
The duality functor $\Dual$ identifies this with the co-
(resp.~contra-) variant submodule $\Ordt$ (resp.~$\Ordi$) of $\Ordit$.

\begin{remark}
  \label{rem:universal-comodule}
  The initial (and zero) object of $\Ordt$ (resp.~$\Ordi$) is
  universal among modules (or ``comodules'', to emphasize the
  categorical variance) over the comonoid $\Bir\Univ$ in a co-
  (resp.~contra-) variant $\Ordit$-module.
  Indeed, it corresponds to the universal co- (resp.~contra-) variant
  module in $\Ordt^\op$ (resp.~$\Ordi^\op$) over the \emph{co}monoid
  $\Univ$ in $\Ord^\op$.
\end{remark}

\begin{corollary}\label{cor:adjoint-hom}
  The forgetful functor $\Ordit\into \Ordt \text{(resp.~$\Ordi$)}$
  lifts to a right adjoint of the functor $i\colon \Ordt\to \Ordit$
  (resp.~$t\colon \Ordi\to \Ordit$) in the bicategory of co-
  (resp.~contra-) variant $\Ordit$-modules.
\end{corollary}
\begin{proof}
  This is a corollary of Theorem~\ref{thm:adjoint-homomorphism}.
\end{proof}

\begin{theorem}\label{prop:commuting-actions}
  The co- (resp.~contra-) variant actions on $\Ordt$ (resp.~$\Ordi$)
  by $\Ord$ and by $\Ordit$ commute, to together give an action by the
  monoid $\Ord\times \Ordit$.
\end{theorem}
\begin{proof}
  For $\Ordt$, we need to prove that the functors
  \[
    \Ordit\times \Ord\times \Ordt\xlongrightarrow{\id\times \Disj}
    \Ordit\times \Ordt\xlongrightarrow{\sum}
    \Ordt
  \]
  where $\sum$ denotes the action of $\Ordit$ on its covariant module,
  and
  \[
    \Ord\times \Ordit\times \Ordt\xlongrightarrow{\id\times \sum}
    \Ord\times \Ordt\xlongrightarrow{\Disj}
    \Ordt,
  \]
  are equal after we identify their sources by flipping the first two
  factors.
  This equality is equivalent to that the restriction of those
  functors along any functor $\cimp{1}\to \Ordt$ be equal.

  Any functor $\cimp{1}\to \Ordt$ lifts along the composite
  \begin{equation}
    \label{eq:homomorphism-composite}
    % https://q.uiver.app/#q=WzAsMyxbMCwwLCJcXE9yZFxcdGltZXMgXFxPcmRpdCJdLFsxLDAsIlxcT3JkXFx0aW1lcyBcXE9yZHQiXSxbMiwwLCJcXE9yZHQiXSxbMSwyLCJcXERpc2oiXSxbMCwxLCIiLDEseyJzdHlsZSI6eyJ0YWlsIjp7Im5hbWUiOiJob29rIiwic2lkZSI6InRvcCJ9fX1dXQ==
    \begin{tikzcd}
      {\Ord\times \Ordit} & {\Ord\times \Ordt} & \Ordt
      \arrow[hook, from=1-1, to=1-2]
      \arrow["\Disj", from=1-2, to=1-3]
    \end{tikzcd}
  \end{equation}
  of homomorphisms of covariant $\Ord$-modules, where the action of
  $\Ord$ on $\Ord\times \Ordit$ (or the other direct product in the
  composition) is on the factor $\Ord$.
  Therefore, it suffices to prove that the restriction of $\sum\colon
  \Ordit\times \Ordt\to \Ordt$ along this
  composite~(\ref{eq:homomorphism-composite}) is a homomorphism of
  covariant $\Ord$-modules.

  As a functor, (\ref{eq:homomorphism-composite}) is equal to the
  composite
  \begin{equation}
    \label{eq:factorized-homomorphism}
    \Ordit\times \Ord\xlongrightarrow{\id\times t}
    \Ordit\times \Ordt\xlongrightarrow{\sum}
    \Ordt
  \end{equation}
  after we identify this source with the source of
  (\ref{eq:homomorphism-composite}) by flipping the factors.
  The restriction of $\sum$ along this
  composite~(\ref{eq:factorized-homomorphism}) appears as the upper
  right composite in the commutative diagram
  % https://q.uiver.app/#q=WzAsNixbMiwxLCJcXE9yZHQiXSxbMCwxLCJcXE9yZGl0XFx0aW1lcyBcXE9yZCJdLFswLDAsIlxcT3JkaXRcXHRpbWVzIFxcT3JkaXRcXHRpbWVzIFxcT3JkIl0sWzEsMCwiXFxPcmRpdFxcdGltZXMgXFxPcmRpdFxcdGltZXMgXFxPcmR0Il0sWzEsMSwiXFxPcmRpdFxcdGltZXMgXFxPcmR0Il0sWzIsMCwiXFxPcmRpdFxcdGltZXMgXFxPcmR0Il0sWzIsMywiXFxpZFxcdGltZXMgdCJdLFsyLDEsIlxcc3VtXFx0aW1lcyBcXGlkIiwyXSxbMyw1LCJcXGlkXFx0aW1lcyBcXHN1bSJdLFszLDQsIlxcc3VtXFx0aW1lcyBcXGlkIiwyXSxbNCwwLCJcXHN1bSIsMl0sWzEsNCwiXFxpZFxcdGltZXMgdCIsMl0sWzUsMCwiXFxzdW0iXV0=
\[\begin{tikzcd}
	{\Ordit\times \Ordit\times \Ord} & {\Ordit\times \Ordit\times \Ordt} & {\Ordit\times \Ordt} \\
	{\Ordit\times \Ord} & {\Ordit\times \Ordt} & \Ordt
	\arrow["{\id\times t}", from=1-1, to=1-2]
	\arrow["{\sum\times \id}"', from=1-1, to=2-1]
	\arrow["{\id\times \sum}", from=1-2, to=1-3]
	\arrow["{\sum\times \id}"', from=1-2, to=2-2]
	\arrow["\sum", from=1-3, to=2-3]
	\arrow["{\id\times t}"', from=2-1, to=2-2]
	\arrow["\sum"', from=2-2, to=2-3]
\end{tikzcd}\]
  in which the left vertical functor and the bottom composite, which
  is identical to (\ref{eq:factorized-homomorphism}), are
  homomorphisms of $\Ord$-modules.
  This proves that the restriction of $\sum\colon
  \Ordit\times \Ordt\to \Ordt$ along the
  composite~(\ref{eq:homomorphism-composite}) is a homomorphism of
  $\Ord$-modules as desired.

  The result for $\Ordi$ follows by symmetry or the symmetric
  argument.
\end{proof}

\begin{theorem}
  The pairing $\Disj\colon \Ordi\times \Ordt\to \Ordit$ is a
  homomorphism of $\Ordit$-bimodules.
\end{theorem}
\begin{proof}
  Corollary~\ref{cor:adjoint-hom} implies that the functor is a
  homomorphism of covariant modules since the composite
  % https://q.uiver.app/#q=WzAsMyxbMCwwLCJcXE9yZGlcXHRpbWVzIFxcT3JkdCJdLFsxLDAsIlxcT3JkaXQiXSxbMiwwLCJcXE9yZHQiXSxbMCwxLCJcXERpc2oiXSxbMSwyLCIiLDAseyJzdHlsZSI6eyJ0YWlsIjp7Im5hbWUiOiJob29rIiwic2lkZSI6InRvcCJ9fX1dXQ==
\[\begin{tikzcd}
	{\Ordi\times \Ordt} & \Ordit & \Ordt
	\arrow["\Disj", from=1-1, to=1-2]
	\arrow[hook, from=1-2, to=1-3]
\end{tikzcd}\]
  is such by the construction of the pairing and
  Theorem~\ref{prop:commuting-actions}.
  By symmetry or the symmetric argument, the functor is also a
  homomorphism of contravariant modules.
\end{proof}

From this, we obtain a refinement of
Corollary~\ref{cor:pairing-adjoint-yuru} as follows.

\begin{corollary}
  \label{cor:pairing-adjoint}
  The pairing $\Disj\colon \Ordi\times \Ordt\to \Ordit$ determines a
  homomorphism $\Ordt\to \Hom_{\Mod^\contra_{\Ordit}}(\Ordi, \Ordit)$
  of covariant ($\Ord\times \Ordit$)-modules, where the target is
  homomorphisms $\Ordi\to \Ordit$ of contravariant $\Ordit$-modules.
  Equivalently, the pairing $\Disj$ determines a homomorphism
  $\Ordi\to \Hom_{\Mod^\co_{\Ordit}}(\Ordt, \Ordit)$ of contravariant
  ($\Ord\times \Ordit$)-modules, where the target is homomorphisms of
  covariant $\Ordit$-modules.
\end{corollary}

\begin{remark}
  Each homomorphism of Corollary~\ref{cor:pairing-adjoint-yuru} is
  the composite of the corresponding homomorphism here with the
  $\Ord$-linear forgetful functor from the target here to the target
  over there.
\end{remark}

\begin{remark}
  \label{rem:universality}
  It follows from the universal property of $\Ordt$
  (Remark~\ref{rem:universal-comodule}) that
  $\Hom_{\Mod^\co_{\Ordit}}(\Ordt, \Ordit)$ is equivalent as a
  category (in fact, as a contravariant $\Ordit$-module) to
  $\Mod^\co_{\Bir\Univ}(\Ordit)$, that of covariant modules in
  $\Ordit$ over the universal comonoid $\Bir\Univ$.
  Alternatively, using Remark~\ref{rem:universal-covariant-module},
  the category
  \begin{align}
    \label{eq:duality-on-hom}
    \Hom_{\Mod^\co_{\Ordit}}(\Ordt, \Ordit)
    &\equivwith \Hom_{\Mod^\co_{\Ord^\op}}(\Ordt^\op, \Ord^\op)\\
    &\equivwith \Hom_{\Mod^\co_{\Ord}}(\Ordt, \Ord)^\op
      \notag
  \end{align}
  is equivalent to $\Mod^\co_\Univ(\Ord)^\op$.

  Similarly, we have equivalence
  \[
    \Hom_{\Mod^\contra_{\Ordit}}(\Ordi, \Ordit)
    \equivwith \Mod^\contra_{\Bir\Univ}(\Ordit)
    \equivwith \Mod^\contra_\Univ(\Ord)^\op
  \]
  of categories.
\end{remark}

The pairing $\Disj$ is \emph{perfect} in the following sense.

\begin{theorem}
  \label{thm:perfect}
  The homomorphisms of Corollary~\ref{cor:pairing-adjoint} are in fact
  isomorphisms.
\end{theorem}

For the proof, we use the following fact.

\begin{lemma}
  The underlying object functor $\Mod^\contra_{\Bir\Univ}(\Ordit)\to
  \Ordit$ is equivalent to $i\colon \Ordt\to \Ordit$ as a functor to
  (i.e., category over) $\Ordit$.
\end{lemma}
\begin{proof}
  The free (or ``cofree'', to emphasize the categorical variance)
  contravariant $\Bir\Univ$-module comonad on $\Ordit$ is equal to the
  comonad associated to the adjuction between $i\colon \Ordt\to
  \Ordit$ and the forgetful functor $\Ordt\hookleftarrow \Ordit$.
  Moreover, this adjunction is comonadic by the Barr--Beck theorem.
\end{proof}

\begin{remark}
  The isomorphism will turn out shortly to be linear over $\Ordit$.
  This results in a concrete interpretation of the linearity of the
  adjunction in Corollary~\ref{cor:adjoint-hom}.
\end{remark}

\begin{proof}[Proof of Theorem]
  To show that the homomorphism $\Ordt\to
  \Hom_{\Mod^\contra_{\Ordit}}(\Ordi, \Ordit)$ is an isomorphism, it
  suffices to prove that it is an equivalence on the underlying
  categories.
  Thus, it suffices to prove that the composite
  \begin{equation}
    \label{eq:pairing-to-modules}
    \Ordt
    \xlongrightarrow{\Disj} \Hom_{\Mod^\contra_{\Ordit}}(\Ordi, \Ordit)
    \longequivto \Mod^\contra_{\Bir\Univ}(\Ordit)
  \end{equation}
  is an equivalence of categories.
  See Remark~\ref{rem:universality}.

  For this, it suffices to prove that the composite of
  (\ref{eq:pairing-to-modules}) followed by the underlying object
  functor $\Mod^\contra_{\Bir\Univ}(\Ordit)\to \Ordit$ is equal to
  $i\colon \Ordt\to \Ordit$.
  Indeed, Lemma and monicity of $i\colon \Ordt\to \Ordit$ will then
  imply that (\ref{eq:pairing-to-modules}) is equal to the isomorphism
  of Lemma.

  To compute it, the composite
  \[
    \Hom_{\Mod^\contra_{\Ordit}}(\Ordi, \Ordit)
    \longequivto \Mod^\contra_{\Bir\Univ}(\Ordit)
    \longto \Ordit
  \]
  is evaluation of the underlying functor of each homomorphism
  $\Ordi\to \Ordit$ at the underlying object of the universal
  contravariant $\Bir\Univ$-module in $\Ordi$, which is the initial
  object.
  Moreover, the pairing $\Disj$ on $\Ordt$ with the initial object of
  $\Ordi$, followed by the forgetful functor $\Ordit\into \Ordt$, is
  equal to the composite
  % https://q.uiver.app/#q=WzAsMyxbMCwwLCJcXE9yZHQiXSxbMSwwLCJcXE9yZGl0Il0sWzIsMCwiXFxPcmR0Il0sWzAsMSwiaSJdLFsxLDIsIiIsMCx7InN0eWxlIjp7InRhaWwiOnsibmFtZSI6Imhvb2siLCJzaWRlIjoidG9wIn19fV1d
\[\begin{tikzcd}
	\Ordt & \Ordit & \Ordt
	\arrow["i", from=1-1, to=1-2]
	\arrow[hook, from=1-2, to=1-3]
\end{tikzcd}\]
  by the construction of the pairing.
  The result follows since the forgetful functor $\Ordit\into \Ordt$
  is monic.

  The result for the other homomorphism follows from symmetry or the
  symmetric argument.
\end{proof}

We record the result of a computation done during the proof.

\begin{example}
  \label{ex:homomorphism-by-pairing}
  $i\colon \Ordt\to \Ordit$ is a zero object of
  $\Hom_{\Mod^\co_{\Ordit}}(\Ordt, \Ordit)$ since it is pairing by
  $\Disj$ with the zero object of $\Ordi$.
  Symmetrically, $t\in \Hom_{\Mod^\contra_{\Ordit}}(\Ordi, \Ordit)$ is
  a zero object.
\end{example}

\begin{notation}
  $\sum\colon \Ordi\times \Ordt\to \Ord$ denotes the pairing
  isomorphic to $\Disj\colon \Ordi^\op\times \Ordt^\op = (\Ordi\times
  \Ordt)^\op\to \Ordit^\op$ through the duality isomorphisms.
\end{notation}

\begin{corollary}
  The pairing $\sum$ determines an isomorphism $\Ordt\equivto
  \Hom_{\Mod^\contra_{\Ord}}(\Ordi, \Ord)$ (resp.~$\Ordi\equivto
  \Hom_{\Mod^\co_{\Ord}}(\Ordt, \Ord)$) of co- (resp.~contra-) variant
  ($\Ord\times \Ordit$)-modules.
\end{corollary}
\begin{proof}
  This follows from Corollary~\ref{cor:pairing-adjoint} and
  Theorem~\ref{thm:perfect}.
\end{proof}

\begin{example}
  The forgetful functor $\text{$\Ordt$ (resp.~$\Ordi$)}\into \Ord$ is
  a zero object of $\Hom_{\Mod^\co_{\Ord}}(\Ordt, \Ord)$
  (resp.~$\Hom_{\Mod^\contra_{\Ord}}(\Ordi, \Ord)$) since it is
  pairing by $\sum$ with the zero object of $\Ordi$ (resp.~$\Ordt$).
\end{example}

\begin{remark}
  This conclusion in $\Hom_{\Mod^\co_{\Ord}}(\Ordt, \Ord)$ is of
  course the dual of Example~\ref{ex:homomorphism-by-pairing} via the
  isomorphism~(\ref{eq:duality-on-hom}) (which has implicitly been
  used in the deduction of Corollary from its dual results obtained
  earlier), and similarly for the conclusion in
  $\Hom_{\Mod^\contra_{\Ord}}(\Ordi, \Ord)$.
\end{remark}

\begin{remark}
  As a result of the correspondence seen in Example, the isomorphism
  of Corollary as that of co- (resp.~contra-) variant $\Ord$-modules
  is over $\Ord$, where $\Hom_{\Mod^\contra_{\Ord}}(\Ordi, \Ord)$
  (resp.~$\Hom_{\Mod^\co_{\Ord}}(\Ordt, \Ord)$) is mapped to $\Ord$ by
  evaluation of the underlying functors at the zero object of the
  source.
  ($\Ordt$ and $\Ordi$ are mapped to $\Ord$ by the respective
  forgetful functors.)
\end{remark}

\begin{remark}
  As in Remark~\ref{rem:universality},
  $\Hom_{\Mod^\contra_{\Ord}}(\Ordi, \Ord)$
  (resp.~$\Hom_{\Mod^\co_{\Ord}}(\Ordt, \Ord)$) as a co-
  (resp.~contra-) variant $\Ord$-modules over $\Ord$ is isomorphic to
  $\Mod^{variance}_\Univ(\Ord)$, where $variance$ is $\contra$
  (resp.~$\co$).
\end{remark}

\begin{remark}
  $\Mod^{variance}_\Univ(\Ord)$ with $variance$ $\contra$
  (resp.~$\co$) is monadic over $\Ord$, where the adjunction/monad is
  co- (resp.~contra-) variantly linear over $\Ord$.
  Together with the previous Remarks, this results in a concrete
  interpretation of the linearity of the adjunctions in
  Theorem~\ref{thm:adjoint-homomorphism}
  (resp.~Corollary~\ref{cor::adjoint-homomorphism}).
\end{remark}


\begin{thebibliography}{99}
\bibitem{modulate}Carboni, Aurelio; Johnson, Scott; Street, Ross; Verity, Dominic. Modulated bicategories, Journal of Pure and Applied Algebra \textbf{94} (1994) 229--282.
\bibitem{enrich}Kelly, Gregory Maxwell. Basic concepts of enriched
  category theory, London Mathematical Society Lecture Note Series,
  \textbf{64}. Cambridge etc.: Cambridge University Press. (1982).\
\bibitem{riron}Matsuoka, T. \emph{Higher theories of algebraic
    structures}. arXiv:1601.00301 (math.CT)
\end{thebibliography}
\end{document}